\documentclass[leqno,final]{siamltex}
 
\usepackage{graphicx} 
\usepackage{amsmath,amstext,amssymb,bm}
\usepackage{leftidx}

\usepackage{xcolor} 
\usepackage{soul} 
\usepackage{tikz}
\usetikzlibrary{shapes,arrows}
\usepackage{mathrsfs}
\usepackage{epstopdf}
\usepackage{color}
\usepackage{multirow}
\usepackage{tabularx}
\usepackage[shortlabels]{enumitem}

\numberwithin{equation}{section}
\newtheorem{remark}{Remark}[section]

\allowdisplaybreaks[4]

\newcommand{\E}{\mathcal{E}}

\renewcommand{\div}{\mbox{\rm div\,}}

\newcommand{\cF}{\mathcal{F}}

\newcommand{\mH}{\mathbb{H}}
\newcommand{\mP}{\mathbb{P}}
\newcommand{\mE}{\mathbb{E}}
\newcommand{\mV}{\mathbb{V}}

\newcommand{\Ome}{\Omega}
\newcommand{\p}{\partial}
\newcommand{\nab}{\nabla}
\newcommand{\vu}{{\bf u}}

\newcommand{\vW}{{\bf W}}
\newcommand{\vg}{{\bf g}}
\newcommand{\vf}{{\bf f}}
\newcommand{\vH}{{\bf H}}

\newcommand{\vL}{{\bf L}}
\newcommand{\vQ}{{\bf Q}}
\newcommand{\vv}{{\bf v}}
\newcommand{\ve}{{\bf e}}
\newcommand{\vE}{{\bf E}}
\newcommand{\vX}{{\bf X}}
\newcommand{\vw}{{\bf w}}
\newcommand{\vA}{{\bf A}}

\newcommand{\pphi}{\pmb{\phi}}

\begin{document}
	
	\title{High moment and pathwise error estimates for fully discrete mixed finite element approximations of stochastic Navier-Stokes equations with additive noise \thanks{This work was partially supported by the NSF grants DMS-1620168 and DMS-2012414.} }
	\markboth{XIAOBING FENG AND LIET VO}{Stochastic Navier-Stokes equations}

	\author{Xiaobing Feng\thanks{Department of Mathematics, The University of Tennessee, Knoxville, TN 37996, U.S.A. (xfeng@math.utk.edu).}
		\and
		Liet Vo\thanks{Department of Mathematics, The University of Tennessee, Knoxville, TN 37996. {\em Current Address:}			
			Department of Mathematics, Statistics and Computer Science, The University of Illinois at Chicago, Chicago, IL 60607, U.S.A. (lietvo@uic.edu).}
	}

	\maketitle
	
	\begin{abstract} 
		 This paper is concerned with high moment and pathwise error estimates for fully discrete mixed finite element approximations of stochastic Navier-Stokes equations with general additive noise.  The implicit Euler-Maruyama scheme and standard  mixed finite element methods are employed respectively for the time and space discretizations. High moment error estimates for both velocity and a time-averaged pressure approximations in strong $L^2$ and energy norms are obtained, pathwise error estimates 
		 are derived by using the Kolmogorov Theorem.  
		 Unlike their deterministic counterparts, the spatial error constants grow in the order  of  $O(k^{-\frac12})$, where $k$ denotes time step size.  
		 Numerical experiments are also provided to validate the error estimates and their sharpness.
		
	\end{abstract}
	
	

	\section{Introduction}\label{sec-1}
	We consider the  following time-dependent stochastic Navier-Stokes equations:
	\begin{subequations}\label{eq1.1}
		\begin{alignat}{2} \label{eq1.1a}
			d\vu &=\bigl[\nu\Delta \vu - \vu\cdot \nab\vu -\nabla p + { \vf}\bigr] dt +  \vg(t)d W(t)  &&\qquad\mbox{a.s. in}\, D_T,\\
			\div \vu &=0 &&\qquad\mbox{a.s. in}\, D_T,\label{eq1.1b}\\
			\vu(0)&= \vu_0 &&\qquad\mbox{a.s. in}\, D,\label{eq1.1d}
		\end{alignat}
	\end{subequations}
	where $D = (0,L)^2 \subset \mathbb{R}^2 \,$ represents a period of the periodic domain in $\mathbb{R}^2$, $\vu$ and $p$ stand for respectively the velocity field and the pressure of the fluid, $\{W(t); t\geq 0\}$ denotes real-valued
	Wiener process, and $\vf$ is a body force function. In addition, $\vg$ is the diffusion coefficient (see section \ref{sub5.2.2} for its precise definition).  Here we seek periodic-in-space solutions $(\vu,p)$ with period $L$, that is, 
	$\vu(t,{\bf x} + L{\bf e}_i) = \vu(t,{\bf x})$ and $p(t,{\bf x}+L{\bf e}_i)=p(t,{\bf x})$ 
	almost surely 
	and for any $(t, {\bf x})\in (0,T)\times \mathbb{R}^d$  and $1\leq i\leq d$, where 
	$\{\bf e_i\}_{i=1}^d$ denotes the canonical basis of $\mathbb{R}^d$.
	
	  Numerical analysis of \eqref{eq1.1} has been studied by several researchers.  In  \cite{BCP12} the authors established stability and convergence of the standard mixed finite element method of \eqref{eq1.1}. Later in \cite{CP2012} the authors proved the rates of convergence in probability for the velocity approximation in the case of multiplicative noise. The main difficulty for establishing a strong convergence for any numerical approximation of \eqref{eq1.1} is the interplay between nonlinearity and  stochasticity of the equations. To compute or estimate quantities of stochastic interests such as the expectation and moments, all the norms must have another layer of the integration which is the main reason why the classical Gronwall inequality argument fails. To overcome this difficulty,  in \cite{CP2012} the authors introduced a sequence of sub-sample spaces which converges to the sample space under the probability measure. The error estimates on these sub-sample spaces are computed with partial expectations and are considered as weak convergence, and some improved error estimates of the same type were recently obtained in \cite{Qui_2022}. 
	  In \cite{BM18, BM_space} the authors showed the strong $L^2$-convergence of the mixed finite element method for \eqref{eq1.1} by estimating the error estimates on the complements of these sub-sample spaces, which then leads to  the strong convergence with a logarithmic rate.  Moreover, the authors  were able to establish in \cite{BM2021}  strong convergence with a polynomial rate  in the case of a divergence-free additive noise. 
	  We note that all the above mentioned error estimates  for \eqref{eq1.1}  are second moment estimates and most  are only for the velocity approximation.  
	  No high moment and pathwise error estimates have been reported in the literature so far.  These missing error estimates are important to know because they provide different 
	  quantities of stochastic interests in practice.
	  
	  The primary goals of this paper is to fill such a void in the case of general additive noise and to develop the analysis techniques for deriving high moment and pathwise error estimates for numerical nonlinear stochastic PDEs in general. It should be noted that
	  the desired high moment and pathwise error estimates will be obtained for both the velocity and pressure approximations of the fully discrete mixed finite element method for \eqref{eq1.1}.  Our main ideas are to obtain the former based on an exponential 
	  stability estimate, which is inspired by a similar idea first introduced in  \cite{BM2021}, and a bootstrap technique,  and to obtain the latter by using the 
	  Kolmogorov Theorem (see Theorem \ref{kolmogorov}).
	
	 The remainder of this paper is organized as follows. In section \ref{sec-2}, we present some preliminaries including the definition 
	 of variational solutions to \eqref{eq1.1} and the assumptions on the diffusion function $\vg$. In section \ref{section5.3},  we introduce  the time discretization for \eqref{eq1.1} in Algorithm 1 and establish some stability estimates for its solution, including an exponential stability estimate which plays a crucial role in our error analysis for the velocity approximation in Theorems \ref{theorem_semi_chapter5}, and \ref{thm_pathwise_semi}, and the error estimates for the pressure approximation in Theorems \ref{semi_pressure_error} and \ref{thm_semi_pathwise_pressure}.  In section \ref{section5.4}, we formulate the fully discrete mixed finite element method in Algorithm 2. The highlight of this section is to establish the desired  high moment and pathwise error estimates for both velocity and pressure approximations. Finally, 
	 we present three numerical experiment results in section \ref{sec-5} to validate the error estimates and their sharpness. 
	 
	 
	\section{Preliminaries}\label{sec-2}
	\subsection{Notations}\label{sec-2.1}
	Standard function and space notation will be adopted in this paper. 
	Let $\vH^1_0(D)$ denote the subspace of $\vH^1(D)$ whose ${\mathbb R}^d$-valued functions have zero trace on $\p D$, and $(\cdot,\cdot):=(\cdot,\cdot)_D$ denote the standard $L^2$-inner product, with induced norm $\Vert \cdot \Vert$. We also denote ${\bf L}^p_{per}(D)$ and ${\bf H}^{k}_{per}(D)$ as the Lebesgue and Sobolev spaces of the functions that are periodic.
	 $C$ denotes a generic constant which is independent of the mesh parameter $h$ and $k$.
	
	Let $(\Omega,\cF, \{\cF_t\},\mP)$ be a filtered probability space with the probability measure $\mP$, the 
	$\sigma$-algebra $\cF$ and the continuous  filtration $\{\cF_t\} \subset \cF$. For a random variable $v$ 
	defined on $(\Omega,\cF, \{\cF_t\},\mP)$,
	${\mathbb E}[v]$ denotes the expected value of $v$. 
	For a vector space $X$ with norm $\|\cdot\|_{X}$,  and $1 \leq p < \infty$, we define the Bochner space
	$\bigl(L^p(\Omega;X); \|v\|_{L^p(\Omega;X)} \bigr)$, where
	$\|v\|_{L^p(\Omega;X)}:=\bigl({\mathbb E} [ \Vert v \Vert_X^p]\bigr)^{\frac1p}$.
	We also define 
	\begin{align*}
		{\mathbb H} := \bigl\{{\bf v}\in  \vL^2_{per}(D) ;\,\div {\bf v}=0 \mbox{ in }D\, \bigr\}\, , \quad 
		{\mathbb V} :=\bigl\{{\bf v}\in  \vH^1_{per}(D) ;\,\div {\bf v}=0 \mbox{ in }D \bigr\}\, .
	\end{align*}
	
	We recall from \cite{Girault_Raviart86} the (orthogonal) Helmholtz projection 
	${\bf P}_{{\mathbb H}}: \vL^2_{per}(D) \rightarrow {\mathbb H}$  and define the Stokes operator ${\bf A} := -{\bf P}_{\mathbb H} \Delta: {\mathbb V} \cap \vH^2_{per}(D) \rightarrow {\mathbb H}$. 
	
	\subsection{Some useful facts and inequalities} In this subsection, we cite some   useful facts and inequalities which will  be used in later sections.
	
	First of all, we recall the Kolmogorov Criteria/Theorem for a pathwise continuity of stochastic processes (cf. \cite{PZ1992}).
	
	\begin{theorem}[Kolmogorov's criteria]\label{kolmogorov}
		Let $\vX(t), t \in [0,T]$, be a stochastic process with values in a separable Banach space $E$ such that, for some positive constant $C>0$, $\alpha > 0, \beta > 0$ and all $t,s \in [0,T]$,
		\begin{align}
			\mE\bigl[\|\vX(t) - \vX(s)\|^{\beta}\bigr] \leq C|t-s|^{1+\alpha}.
		\end{align}
		Then for each $T>0$, almost every $\omega$ and each $0< \gamma < \frac{\alpha}{\beta}$ there exists a constant $K = K(\omega,\gamma,T)$ such that
		\begin{align}
			\|\vX(t,\omega) - \vX(s,\omega)\| \leq K|t-s|^{\gamma}\qquad\mbox{ for all } t,s \in [0,T].
		\end{align}
		Moreover, $\mE\bigl[|K|^\beta\bigr] < \infty$ for all $\beta >0$.
	\end{theorem}

	\begin{lemma}[Burkholder-Davis-Gundy inequality, \cite{PZ1992}]
		Let $\pphi(t)$ be a stochastic process for all $t \in [0,T]$. For any $p >0$, there exists a positive constant $C_b$ such that: 
		\begin{align}
			\mE\biggl[\max_{0\leq t \leq T}\biggl\|\int_0^t \pphi(\xi)\, dW(\xi)\biggr\|_2^p\biggr] \leq C_b\, \mE\biggl[\biggl(\int_0^T \|\pphi(\xi)\|^2_{L^2}\, d\xi\biggr)^{\frac{p}{2}}\biggr],
		\end{align}
		where $C_b = C_b(T,p)$.
\end{lemma}
	
	The next lemma recalls the well-known It\^o isometry and also states a helpful inequality for stochastic processes 	(cf. \cite{PZ1992}). 
	
	\begin{lemma}\label{lemma2.3}
		Let $\pphi(t)$ be a stochastic process for all $t \in [0,T]$. Define $\displaystyle \vX_t = \int_0^t \pphi(\xi)\, dW(\xi)$. We have:
		\begin{enumerate}[{\rm (i)}]
			\item If $\pphi \in L^2(\Ome;L^2(0,T; \vL^2(D)))$, then
			\begin{align}\label{ito2}
				\mE\bigl[\|\vX_t\|^2_{\vL^2}\bigr] = \mE\biggl[\int_0^t\|\pphi(\xi)\|^2_{\vL^2}\, d\xi\biggr].
			\end{align}
			\item If $\pphi \in L^p(\Ome;L^p(0,T;\vL^2(D)))$ for $p > 2$, then
			\begin{align}\label{ito4}
				\mE\bigl[\|\vX_t\|^p_{\vL^2}\bigr] \leq C(t,p)\,\mE\biggl[\int_0^t \|\pphi(\xi)\|^p_{\vL^2}\, d\xi\biggr],
			\end{align}
			where $\displaystyle C(t,p) = \frac{C_b}{2}(p-1)(p-2) t^{\frac{p}{2}} + (p-1)C_b$.
		\end{enumerate}
	\end{lemma}

The next lemma states  both continuous and discrete Gronwall inequalities.

	\begin{lemma}[Gronwall inequalities]\label{gronwall}
		\begin{enumerate}[{\rm (a)}]
			\item Let $y$ and $g$ be nonnegative and integrable and $c$ a nonnegative constant. If 
			\begin{align*}
				y(t) \leq c + \int_0^t g(s)y(s)\, ds\qquad\forall t \geq 0,
			\end{align*}
			then $\displaystyle y(t) \leq c\exp\biggl(\int_0^t g(s)\, ds\biggr)$ for all $t \geq 0$.
			
			\item Let $\{y_n\}$ and $\{g_n\}$ be nonnegative sequences and $c$ a nonnegative constant. If 
			\begin{align*}
				y_n \leq c + \sum_{0\leq j < n} g_jy_j\qquad\forall n \geq 0,
			\end{align*}
			then $\displaystyle y_n \leq c \prod_{0\leq j <n}(1+g_j) \leq c\exp\biggl(\sum_{0\leq j < n}g_j\biggr)$ for all $n \geq 0$.
		\end{enumerate}
	\end{lemma}	
	
	\smallskip
	
	\begin{lemma}[Hoeffding's lemma]
		Let $X$ be a real-valued random variable with $\mE[X] = \mu$, such that $a\leq X \leq b$ a.s., for some $a,b \in \mathbb{R}$. Then, for all $\lambda \in \mathbb{R}$,
		\begin{align}
			\mE\bigl[e^{\lambda X}\bigr] \leq \exp\Bigl(\lambda\mu + \frac{\lambda^2}{8}(b-a)^2\Bigr).
		\end{align}
	\end{lemma}
	
	Finally, we recall the following property of  the $\mathbb{R}$-valued Wiener process.
	 
	\begin{align}\label{mean_wiener}
		\mE\Bigl[|W(t) - W(s)|^{2m}\Bigr] \leq C_m|t-s|^{m}\qquad\forall m \in \mathbb{N},
	\end{align}
	where, for $m=1$, we obtain the equality and $C_m =1$. See \cite{Ichikawa} for higher dimensional Wiener process.
	
 In addition, there exists a deterministic constant $C>0$ such that:
	\begin{align}\label{a.s_wiener}
		|W(t) - W(s)| \leq C|t-s|^{\alpha},\qquad\mP-a.s.,
	\end{align}
	where,  when $m=1$,  the equality holds and $C_1 =1$. See \cite{Ichikawa} for high dimensional generalizations.
	 
	\subsection{Variational solutions}\label{sub5.2.2}
	Let $\vg(\cdot,\cdot): D\times [0,T] \rightarrow \mathbb{R}$ and $\vg \in 
	L^{\infty}\bigl(0,T;\vH^1(D)\bigr)$. We make the following additional assumptions on $\vg$:
	\begin{enumerate}[(G1)]
		\item $\vg$ is Lipschitz continuous in time, that is, there is a constants $C_{\vg} > 0$ such that
		\begin{align*}
			\|\vg(t) - \vg(s)\|_{\vL^2} \leq C_{\vg}|t-s|\qquad\forall t,s \in [0,T].
		\end{align*}
		\item Let $\sigma >0$. Assume that there exists a constant $ 0< K_0 \leq \frac{\nu}{8 C_L\sigma}$ such that $\|\vg\|_{L^{\infty}(0,T;\vH^1)} \leq K_0$, where ${C}_L$ is the constant from the inequality \eqref{interpolation}.
	\end{enumerate}

\smallskip
	\begin{remark}
		The assumption (G2) will be crucially used to obtain the exponential stability estimates in Lemma \ref{lemma_exp_discrete}. On the other hand, both (G1) and (G2) will be needed to derive the order of convergence in Theorem \ref{theorem_semi_chapter5} and Theorem \ref{theorem_fully_chapter5}. 
	\end{remark}
	
	Next, we introduce the variational solution concept for \eqref{eq1.1} and refer the reader to \cite{Chow07,PZ1992} for a proof of its existence and uniqueness.
	
	\smallskip
	\begin{definition}\label{def2.1} 
		Given $(\Omega,\cF, \{\cF_t\},\mP)$, let $W$ be an ${\mathbb R}$-valued Wiener process on it. 
		Suppose ${\bf u}_0\in L^2(\Omega, {\mathbb V})$.
		An $\{\cF_t\}$-adapted stochastic process  $\{{\bf u}(t) ; 0\leq t\leq T\}$ is called
		a variational solution of \eqref{eq1.1} if ${\bf u} \in  L^2\bigl(\Omega; C([0,T]; {\mathbb V})) 
		\cap L^2\bigl(\Ome;0,T;\vH^2_{per}(D)\bigr)$,
		and satisfies $\mP$-a.s.~for all $t\in (0,T]$
		\begin{align}\label{equu2.8a}
			\bigl({\bf u}(t),  {\bf v} \bigr)& + \int_0^t  \nu\bigl(\nab {\bf u}(s), \nab {\bf v} \bigr) 
			\,  ds + \int_0^t \big(\vu(s)\cdot\nab\vu(s),\vv\big)\, ds
			\\\nonumber
			&=({\bf u}_0, {\bf v}) + {  \Bigl(\int_0^t \vg(s)\, dW(s), {\bf v} \Bigr)}  \qquad\forall  \, {\bf v}\in {\mathbb V}\, . 
		\end{align}
	\end{definition}
	
	We also recall some properties of the convection term. Define
	\begin{align*}
		b(\vu, \vv, \vw) := \bigl(\vu\cdot\nab\vv,\vw\bigr)\qquad \forall \vu, \vv,\vw \in \vH^1_{per}(D).
	\end{align*}
	It is easy to check that $b(\vu,\cdot,\cdot)$ is skew-symmetric for $\vu \in \mV$, that is,
	\begin{align}
		b(\vu,\vv,\vw) = - b(\vu,\vw,\vv)\qquad\forall \vv, \vw \in \vH^1_{per}
	\end{align}
	and 
	\begin{align*}
		b(\vu,\vv,\vv) = 0 \qquad\forall \vv \in \vH^1_{per}(D).
	\end{align*}
	
	In addition, we also recall the following interpolation inequality for $\vH^1$-functions:
	\begin{align}\label{interpolation}
		b(\vu,\vv,\vw) 
		\leq C_L\|\vu\|^{1/2}_{\vL^2}\|\nab\vu\|^{1/2}_{\vL^2}\|\nab\vv\|_{\vL^2} \|\vw\|^{1/2}_{\vL^2}\|\nab\vw\|^{1/2}_{\vL^2}.
	\end{align}
	
	Definition \ref{def2.1} only defines the velocity $\mathbf{u}$ for \eqref{eq1.1}, 
	its associated pressure $p$ is subtle to define. We do so in the following theorem.
	
	\begin{theorem}\label{thm 2.2}
		Let $\{{\bf u}(t) ; 0\leq t\leq T\}$ be a variational solution of \eqref{eq1.1}. There exists a unique adapted process 
	 $P\in {L^2\bigl(\Omega; L^2(0,T; H^1_{per}(D)/\mathbb{R})\bigr)}$ such that $(\mathbf{u}, P)$ satisfies 
		$\mP$-a.s.~for all $t\in (0,T]$
		\begin{subequations}\label{equu2.100}
			\begin{align}\label{equu2.10a}
				&\bigl({\bf u}(t),  {\bf v} \bigr) + \nu\int_0^t  \bigl(\nab {\bf u}(s), \nab {\bf v} \bigr) \, ds + \int_0^t \big(\vu(s)\cdot\nab\vu(s),\vv\big)\, ds
				- \bigl(  \div \mathbf{v}, P(t) \bigr) \\
				&=({\bf u}_0, {\bf v}) + \int_0^t \big(\vf(s), \vv\big) \, ds 
				+  {\int_0^t  \bigl( {\bf B}\bigl({\bf u}(s)\bigr), {\bf v} \bigr)\, dW(s)}  \,\,\, \forall  \, {\bf v}\in \vH^1_{per}(D)\, , \nonumber \\ 
				&\bigl(\div {\bf u}, q \bigr) =0 \qquad\forall \, q\in  L^2_{per}(D)/\mathbb{R} .  \label{equu2.10b}
			\end{align}
		\end{subequations}
	\end{theorem}
Since the proof is similar \cite[Theorem 1.3]{FPL2021}, so we omit it.

	\begin{lemma}[\cite{Breit,CP2012}]\label{stability_pdes} 
		Let $\vu$ be a solution defined in Definition \ref{def2.1}. Then we have
		\begin{enumerate}[{\rm (a)}]
			\item Assume that $\vu_0 \in L^r\bigl(\Ome; \mV\bigr)$ for some $r \geq 2$. Then there hold 
			\begin{align*}
				\mE\biggl[\Bigl(\sup_{0\leq t \leq T} \|\nab\vu(t)\|^2_{\vL^2} + \int_0^T \nu\|\nab^2\vu(t)\|^2_{\vL^2}\, dt\Bigr)^{\frac{r}{2}} \biggr] \leq C_r \mE\Bigl[\|\nab \vu_0\|^r_{\vL^2}\Bigr].
			\end{align*}
			\item Assume that $\vu_0 \in L^r\bigl(\Ome; \mV\cap\vH^2(D)\bigr)\cap L^{5r}\bigl(\Ome,\mV\bigr)$ for some $r \geq 2$. Then we have
			\begin{align*}
				\mE\biggl[\Bigl(\sup_{0\leq t \leq T} \|\nab^2\vu(t)\|^2_{\vL^2} + \int_0^T \nu\|\nab^3\vu(t)\|^2_{\vL^2}\, dt\Bigr)^{\frac{r}{2}} \biggr] \leq C_r \mE\biggl[\Bigl(\|\vu_0\|^2_{\vH^2} + \|\nab\vu_0\|^{10}_{\vL^2}\Bigr)^{\frac{r}{2}}\biggr].
			\end{align*}
		\end{enumerate}
	\end{lemma}
	
	\subsection{H\"older continuity of the variational solution}
	We cite the following high moment H\"older continuity estimate for the variational solution whose proof
	can be found in \cite{Breit,CHP2012}. 
	\begin{lemma}[\cite{Breit,CHP2012}]\label{lemma2.2}
		Suppose ${\bf u}_0$ satisfies the assumptions in Lemma \ref{stability_pdes} (b) for some $r >2$ and $\vg \in L^2(\Ome;L^{\infty}(0,T;\vH_{per}^2(D)))$. Then there exists a constant $C \equiv C(D_T, \vu_0, \vg,r)>0$, such that the variational solution to problem \eqref{eq1.1} satisfies
		for $s,t \in [0,T]$
		\begin{align}\label{equu2.20a}
			\mE\Bigl[\|\vu(t) - \vu(s)\|^{r}_{\mV}\Bigr] \leq C|t-s|^{r\gamma}\qquad\forall \gamma \in \Bigl(0,\frac12\Bigr).
		\end{align}
	\end{lemma}

	
	\section{Semi-discretization in time}\label{section5.3} 
	In this section, we consider the time discretization of \eqref{equu2.8a} that is based on the fully implicit Euler-Maruyama method.
	
	\subsection{Formulation and stability of the time discretization}\,
	The Euler-Maruyama time discretization of \eqref{equu2.8a} is given by the following algorithm.
 
	\textbf{Algorithm 1.} Let $\vu^0 = \vu_0$ be a given $\mV$-valued random variable. Find $\bigl(\vu^{n+1}, p^{n+1}\bigr) \in L^2(\Ome;\mV \times L^2_{per}(D)/\mathbb{R})$ such that $\mP$-a.s.
	\begin{align}\label{equu3.1}
		\bigl(\vu^{n+1} - \vu^{n},\pphi\bigr) + \nu k\,\bigl(\nab\vu^{n+1},\nab\pphi\bigr) + k\,\bigl(\vu^{n+1}\cdot\nab\vu^{n+1},\pphi\bigr)
		&- k\, \bigl(p^{n+1},\div \pphi\bigr) \\\nonumber
		&= \bigl(\vg^n\Delta W_{n+1},\pphi\bigr),\\
		\bigl(\div \vu^{n+1},\psi\bigr) &=0
	\end{align}
	$\forall \pphi \in \vH^1_{per}(D)$ and $\psi \in L^2_{per}(D)$. Where $\vg^n = \vg(t_n)$ and $\Delta W_{n+1} = W(t_{n+1})- W(t_n)$.
	
	If $\pphi \in \mV$, then \eqref{equu3.1} reduces to
	\begin{align}\label{eq_reforms}
		\bigl(\vu^{n+1} - \vu^{n},\pphi\bigr) + \nu k\,\bigl(\nab\vu^{n+1},\nab\pphi\bigr) + k\,\bigl(\vu^{n+1}\cdot\nab\vu^{n+1},\pphi\bigr) 
		= \bigl(\vg^n\Delta W_{n+1},\pphi\bigr).
	\end{align}
	
	We quote the following stability estimates for $\{\vu^n\}$ given in \cite{BCP12}.
	
	\begin{lemma}\label{stability_means}
		Let $\vu_0 \in L^{2^q}(\Ome;\mV)$ for an integer $1 \leq q <\infty$ be given, such that $\mE\bigl[\|\vu_0\|^{2^q}_{\mV}\bigr] \leq C$. Then there exists a constant $C_{T,q} = C(T, q, \vu_0)$ such that the following estimations hold:
		\begin{enumerate}[{\rm (i)}]
			\item $\displaystyle\mE\biggl[\max_{1\leq n \leq M}\|\vu^n\|^{2^q}_{\mV} + \nu k\sum_{n=1}^M \|\vu^n\|^{2^q-2}_{\mV}\|{\bf A}\vu^n\|^2_{\vL^2}\biggr] \leq C_{T,q}$.
			\item $\displaystyle\mE\biggl[\sum_{n=1}^M \|\vu^n - \vu^{n-1}\|^2_{\mV}\|\vu^n\|^2_{\mV}\biggr] \leq C_{T,2}$.
			\item $\displaystyle\mE\Biggl[\biggl(\sum_{n=1}^M \|\vu^{n} - \vu^{n-1}\|^2_{\mV}\biggr)^q + \biggl(\nu k\sum_{n=1}^M \|{\bf A}\vu^{n}\|^2_{\mV}\biggr)^q\Biggr] \leq C_{T,q}$.
		\end{enumerate}
	\end{lemma}
	
	Stability estimates for the pressure approximation $\{p^n\}$ can also be obtained accordingly using  the stochastic {\em inf-sup} 
	estimate.
	
	\begin{lemma}\label{stability_pressures}
		Let $\{p^n\}_{n=1}^M$ be the pressure approximation from Algorithm 1. Under the assumptions of Lemma \ref{stability_means} and for $1\leq q < \infty$, we have
		\begin{align}\label{stability_pressure2}
			\mE\biggl[\biggl(k\sum_{n=1}^M \|\nab p^{n}\|^2_{\vL^2}\biggr)^q\biggr] \leq \frac{C_{T,q}}{k^q}.
		\end{align} 
	\end{lemma}
	
	Next, we present an exponential stability estimate for the velocity approximation from Algorithm 1. Such an estimate was first 
	derived  in \cite[Theorem 8.3]{BM_space}. However, here we provide a new proof which eliminates the restriction that $M$ must  be sufficiently large required in \cite{BM_space}. This estimate will be used in the proofs of Theorem \ref{theorem_semi_chapter5} and Theorem \ref{theorem_fully_chapter5}.
	
	\begin{lemma}\label{lemma_exp_discrete}
		Assume that $\mE[\exp(4\sigma \|\nab\vu_0\|^2_{\vL^2})] \leq C$ for some $\sigma > 0$. Let $\{\vu^n\}$ be generated by Algorithm 1 and assume that $\vg$ satisfies $(G2)$. Then, there holds
		\begin{align}
			\mE\Bigl[\exp\Bigl(&\sigma\max_{1\leq \ell \leq M}\|\nab\vu^{\ell}\|^2_{\vL^2} + \sigma\nu k\sum_{n=1}^M\|\vA\vu^n\|^2_{\vL^2}\Bigr)\Bigr] \leq C_2,
		\end{align}
		where $C_2 = C_2(\sigma,\nu,T,\mE[\exp(4\sigma \|\nab\vu_0\|^2_{\vL^2})])$.
	\end{lemma}

 \begin{proof}
		Choosing $\pphi = \vA\vu^{n+1}$ in \eqref{eq_reforms} and using integration by parts, the binomial formula $2(a,a-b) = \|a\|^2 - \|b\|^2 + \|a-b\|^2$, we obtain
		\begin{align}\label{equu3.13}
			\|\nab\vu^{n+1}\|^2_{\vL^2} &- \|\nab\vu^n\|^2_{\vL^2} + \|\nab(\vu^{n+1} - \vu^n)\|^2_{\vL^2} + 2\nu k \|\vA\vu^{n+1}\|^2_{\vL^2} \\\nonumber
			&=2\bigl(\nab\vg^{n}\Delta W_{n+1}, \nab(\vu^{n+1}-\vu^n)\bigr) +2\bigl(\nab\vg^{n}\Delta W_{n+1}, \nab\vu^n\bigr)\\\nonumber
			&\leq \|\nab\vg^n\Delta W_{n+1}\|^2_{\vL^2} + \|\nab(\vu^{n+1} - \vu^n)\|^2_{\vL^2} \\\nonumber
			&\qquad\qquad\qquad\qquad+ 2\bigl(\nab\vg^{n}\Delta W_{n+1}, \nab\vu^n\bigr).
		\end{align}
		Therefore,
		\begin{align}
			\|\nab\vu^{n+1}\|^2_{\vL^2} &- \|\nab\vu^n\|^2_{\vL^2} + 2\nu k \|\vA\vu^{n+1}\|^2_{\vL^2} \\\nonumber
			&\leq \|\nab\vg^n\Delta W_{n+1}\|^2_{\vL^2} + 2\bigl(\nab\vg^{n}\Delta W_{n+1}, \nab\vu^n\bigr).
		\end{align}
		Next, we lower the index by one and apply the summation $ \sum_{n=1}^{\ell}$ for any $1 \leq \ell \leq M$ to get
		\begin{align}\label{equu3.18}
			\|\nab\vu^{\ell}\|^2_{\vL^2} +2\nu k\sum_{n=1}^{\ell} \|\vA\vu^{n}\|^2_{\vL^2} &\leq \|\nab\vu^0\|^2_{\vL^2} + \sum_{n=1}^{\ell}\|\nab\vg^{n-1}\Delta W_{n}\|^2_{\vL^2}\\\nonumber
			&+ 2\sum_{n=1}^{\ell}\bigl(\nab\vg^{n}\Delta W_{n+1}, \nab\vu^n\bigr).
		\end{align}
		
		Taking $\displaystyle \max_{1\leq \ell \leq M}$ on both sides of \eqref{equu3.18} leads to 
		\begin{align}\label{equu3.19}
			\max_{1\leq \ell \leq M}\|\nab\vu^{\ell}\|^2_{\vL^2} +2\nu k\sum_{n=1}^{M} \|\vA\vu^{n}\|^2_{\vL^2} &\leq \|\nab\vu^0\|^2_{\vL^2} + \sum_{n=1}^{M}\|\nab\vg^{n-1}\Delta W_n\|^2_{\vL^2}\\\nonumber
			&+2\max_{1\leq \ell \leq M}\sum_{n=1}^{\ell}\bigl(\nab\vg^{n-1}\Delta W_n,\nab\vu^{n-1}\bigr)\\\nonumber
			&\leq 2\|\nab\vu^0\|^2_{\vL^2} + 2K_0 \sum_{n=1}^{M}|\Delta W_n|^2\\\nonumber
			&+2\max_{2\leq \ell \leq M}\sum_{n=2}^{\ell}\bigl(\nab\vg^{n-1}\Delta W_n,\nab\vu^{n-1}\bigr).
		\end{align}
		For any $\sigma > 0$, denote $Z_{\ell} =2\sigma\sum_{n=2}^{\ell}\bigl(\nab\vg^{n-1}\Delta W_n,\nab\vu^{n-1}\bigr)$, for $2\leq \ell\leq M$. Define the following martingale by using a piece-wise linear function that defines as follow: for $s \in [t_n,t_{n+1})$, with $n =1,2,\cdots, M-1$, set $\underline{s} = t_n$, $\vu^{\underline{s}} = \vu^n$. Then $\displaystyle Z_{\ell} = \tilde{Z}_{t_{\ell}}$, where
		\begin{align}
			\tilde{Z}_{t} = 2\sigma\int_{t_1}^{t}\bigl(\nab\vg^{\underline{s}}\,d W(s),\nab \vu^{\underline{s}}\bigr)\qquad\forall t\in [t_1,T].
		\end{align}
		Then, $\tilde{Z}_t$ is a $\mathcal{F}_t$-martingale for $t \in [t_1,T]$ and its quadratic variation satisfies
		\begin{align}\label{equu3.21}
			\bigl<\tilde{Z}\bigr>_{t_{\ell}}  &\leq 4\sigma^2\int_{t_1}^{t_{\ell}} \|\nab\vg^{\underline{s}}\|^2_{\vL^2}\|\nab\vu^{\underline{s}}\|^2_{\vL^2}\, ds\\\nonumber
			&\leq 4C_L\sigma^2K_0\, k\sum_{n=1}^{\ell}\|\vA\vu^{n}\|^2_{\vL^2},
		\end{align}
		where $C_L$ is the constant from the inequality $\|\nab\vu^n\|^2_{\vL^2} \leq C_L\|\vA\vu^n\|^2_{\vL^2}$.
		
		We multiply \eqref{equu3.19} by $\sigma >0$ and then take the exponential on both sides to get
		\begin{align}\label{equu323}
			\exp\Bigl(&\sigma\max_{1\leq \ell \leq M}\|\nab\vu^{\ell}\|^2_{\vL^2} + \sigma\nu k\sum_{n=1}^M\|\vA\vu^n\|^2_{\vL^2}\Bigr) \\\nonumber
			&\leq \exp\bigl(2\sigma\|\nab\vu^0\|^2_{\vL^2}\bigr)\times \exp\Bigl(2K_0\sigma\sum_{n=1}^M|\Delta W_n|^2\Bigr)\\\nonumber
			&\times\exp\Bigl(\max_{2\leq \ell \leq M}\bigl[Z_{\ell} -2\bigl<\tilde{Z}\bigr>_{t_{\ell}}\bigr]\Bigr)\\\nonumber
			&\times \exp\biggl(2\max_{2\leq \ell \leq M}\bigl<\tilde{Z}\bigr>_{t_{\ell}} - \sigma\nu k\sum_{n=1}^M\|\vA\vu^n\|^2_{\vL^2}\biggr).
		\end{align} 
		Moreover, from \eqref{equu3.21} we have
		\begin{align}\label{equu3.24}
			&2\max_{2\leq \ell \leq M}\bigl<\tilde{Z}\bigr>_{t_{\ell}} - \sigma\nu k\sum_{n=1}^M\|\vA\vu^n\|^2_{\vL^2} \\\nonumber
			&\leq 8 C_L\sigma^2K_0\, k\sum_{n=1}^{\ell}\|\vA\vu^{n}\|^2_{\vL^2} - \sigma\nu k\sum_{n=1}^M\|\vA\vu^n\|^2_{\vL^2}\\\nonumber
			&\leq \bigl(8 K_0C_L\sigma^2 - \nu\sigma\bigr) k\sum_{n=1}^M\|\vA\vu^n\|^2_{\vL^2}.
		\end{align}	
		Now, by using (G2), we have $8 K_0C_L\sigma^2 - \nu\sigma \leq 0$,
		which also implies that 
		\begin{align}\label{equu3.26}
			\exp\biggl(2\max_{2\leq \ell \leq M}\bigl<\tilde{Z}\bigr>_{t_{\ell}} - \sigma\nu k\sum_{n=1}^M\|\vA\vu^n\|^2_{\vL^2}\biggr)  \leq 1\qquad\mP-a.s.
		\end{align}
		Next, taking the expectation on \eqref{equu323} and using \eqref{equu3.26} we obtain
		\begin{align}\label{equu3.27}
			\mE\Bigl[\exp\Bigl(&\sigma\max_{1\leq \ell \leq M}\|\nab\vu^{\ell}\|^2_{\vL^2} + \sigma\nu k\sum_{n=1}^M\|\vA\vu^n\|^2_{\vL^2}\Bigr)\Bigr] \\\nonumber
			&\leq \mE\Bigl[\exp\bigl(2\sigma\|\nab\vu^0\|^2_{\vL^2}\bigr)\times \exp\Bigl(2K_0\sigma\sum_{n=1}^M|\Delta W_n|^2\Bigr)\\\nonumber
			&\times\exp\Bigl(\max_{2\leq \ell \leq M}\bigl[\tilde{Z}_{t_{\ell}} -2\bigl<\tilde{Z}\bigr>_{t_{\ell}}\bigr]\Bigr)\Bigr].
		\end{align}
		Using the H\"older inequality to separate the product on the right-hand side of \eqref{equu3.27} yields
		\begin{align}\label{equu3.28}
			\mE\Bigl[\exp\Bigl(&\sigma\max_{1\leq \ell \leq M}\|\nab\vu^{\ell}\|^2_{\vL^2} + \sigma\nu k\sum_{n=1}^M\|\vA\vu^n\|^2_{\vL^2}\Bigr)\Bigr]\\\nonumber
			&\leq \Bigl(\mE\bigl[\exp\bigl(4\sigma\|\nab\vu^0\|^2_{\vL^2}\bigr)\bigr]\Bigr)^{1/2}\\\nonumber
			&\qquad\times \biggl(\mE\Bigl[\exp\Bigl(8K_0\sigma\sum_{n=1}^M|\Delta W_n|^2\Bigr)\Bigr]\biggr)^{1/4}\\\nonumber
			&\qquad\times \Bigl(\mE\Bigl[\exp\Bigl(\max_{2\leq \ell \leq M}\bigl[4\tilde{Z}_{t_{\ell}} -\frac{1}{2}\bigl<4\tilde{Z}\bigr>_{t_{\ell}}\bigr]\Bigr)\Bigr]\Bigr)^{1/4}.
		\end{align}
		Since $\Bigl\{\exp\Bigl(4\tilde{Z}_{t} -\frac{1}{2}\bigl<4\tilde{Z}\bigr>_{t}\Bigr)\Bigr\}$ is an exponential martingale on $[t_1,T]$, so it suffices to show that $\mE\Bigl[\exp\Bigl(8K_0\sigma\sum_{n=1}^M|\Delta W_n|^2\Bigr)\Bigr]$ is bounded. To the end, let $\beta = 8K_0\sigma$ and $X_n = |\Delta W_n|^2$,\,$ 1 \leq n \leq M$. 
		Then, by using the independence of the increments of $W(t)$, we can rewrite
		\begin{align}\label{equation3.23}
			\mE\Bigl[\exp\Bigl(8K_0\sigma\sum_{n=1}^M|\Delta W_n|^2\Bigr)\Bigr] &= \mE\Bigl[\exp\Bigl(\beta \sum_{n=1}^M X_n\Bigr)\Bigr] = \prod_{n=1}^M \mE\bigl[\exp(\beta  X_n)\bigr].
		\end{align}
		In addition, by using \eqref{mean_wiener}, we have that for each $1\leq n\leq M$, $\mE[X_n] = \mE[|\Delta W_n|^2] = k$ and then by \eqref{a.s_wiener}, we also have $0 \leq X_n \leq C k^{\alpha}$\,\, a.s. \, $\forall \alpha \in (0,1)$. 
		
		Next, we use Hoeffding's lemma to conclude that
		\begin{align}\label{equation3.24}
			\mE\bigl[\exp(\beta X_n)\bigr] \leq \exp\bigl(\beta  \mE[X_n] + \frac{C^2\beta^2}{8} k ^{2\alpha}\bigr) = \exp\bigl(\beta k + \frac{C^2\beta^2}{8} k ^{2\alpha}\bigr).
		\end{align}
		Substituting \eqref{equation3.24} into the right-hand side of \eqref{equation3.23}, we arrive at
		\begin{align*}
			\mE\Bigl[\exp\Bigl(8K_0\sigma\sum_{n=1}^M|\Delta W_n|^2\Bigr)\Bigr] \leq  \exp\bigl(\beta M k + \frac{C^2\beta^2}{8} M k ^{2\alpha}\bigr)  =  \exp\bigl(\beta T + \frac{C^2\beta^2}{8} T k ^{2\alpha-1}\bigr).
		\end{align*}
		
		The proof is complete by choosing $\frac12 \leq \alpha < 1$.
	\end{proof}
	
	\subsection{High moment and pathwise error estimates for the velocity approximation}
	In this subsection, we present the first main result of this paper which establishes the optimal order 
	high moment and sub-optimal order pathwise error estimates for the velocity approximation generated by Algorithm 1.
	
	\begin{theorem}\label{theorem_semi_chapter5} 
		Let $\vu$ be the variational solution to \eqref{equu2.100} and $\{\vu^{n}\}_{n=1}^M$ be generated by  Algorithm 1. Assume that $\vu_0 \in L^{q}(\Ome; \mV)$ and . Then there exists $C_q = C(T,q,\vu_0,\vf)>0$ for any {integer} $2 \leq q < \infty$ and {real number} $0 < \gamma < \frac12$ such that
		\begin{align}\label{equu310}
			\bigl(\mE\bigl[\max_{1\leq n\leq M}\|\vu(t_n) - \vu^n\|^{q}_{\vL^2}\bigr]\bigr)^{\frac1q} + \biggl(\mE\biggl[\Bigl(\nu k \sum_{n = 1}^M \|\nab(\vu(t_n) &- \vu^n)\|^2_{\vL^2}\Bigr)^{q/2}\biggr]\biggr)^{\frac1q} \\\nonumber
			&\leq C_q\, k^{\frac{1}{2}-\gamma}.
		\end{align}
	\end{theorem}
	
	\smallskip
	
  The proof of Theorem \ref{theorem_semi_chapter5} follows immediately  from Lemma \ref{lemma3.4} and \ref{lemma3.5} below.
	
	\smallskip
	
	\begin{lemma}\label{lemma3.4}
		Let $\ve^n := \vu(t_n) - \vu^n$. 	Under the assumptions of Theorem \ref{theorem_semi_chapter5}, there holds $\mP$-a.s.
		\begin{align}
			\max_{1\leq \ell \leq M}\|\ve^{\ell}\|^2_{\vL^2} &+ {\nu k}\sum_{n=1}^{M} \|\nab\ve^{n}\|^2_{\vL^2} \leq \bigl(A + B\bigr) \exp\biggl(k\sum_{n=1}^{M}D^n\biggr),
		\end{align}
		where
		\begin{align*}\label{equu3.288}
			A &:= 2\nu\sum_{n=1}^{M}\int_{t_{n-1}}^{t_{n}}\|\nab(\vu(t_{n}) - \vu(s))\|^2_{\vL^2}\, ds \\\nonumber
			&\qquad+ \frac{4}{\nu}\sum_{n=1}^{M}\int_{t_{n-1}}^{t_{n}} C_L^2\, \|\vu(s) - \vu(t_{n})\|^2_{\vH^1}\|\nab\vu(s)\|^2_{\vL^2} \, ds \\\nonumber
			&\qquad+ \frac{4}{\nu}\sum_{n=1}^{M}\int_{t_{n-1}}^{t_{n}} C_L^2\, \|\vu(s) - \vu(t_{n})\|^2_{\vH^1}\|\nab\vu(t_{n})\|^2_{\vL^2} \, ds\\\nonumber
			&\qquad+ \frac{2C_L^2}{\nu} k \max_{1\leq \ell \leq M}\|\ve^{\ell}\|^2_{\vL^2}\|\nab\vu^{\ell}\|^2_{\vL^2}, 
		\end{align*}
		\begin{align*}
			B &:= 2\sum_{n=1}^{M}\Bigl\|\int_{t_{n-1}}^{t_{n}}\bigl(\vg(s) - \vg(t_{n-1})\bigr)\, dW(s)\Bigr\|_{\vL^2}^2\\\nonumber
			&\qquad+ 2\max_{1\leq \ell \leq M}\Bigl|\sum_{n=1}^{\ell}\biggl(\int_{t_{n-1}}^{t_{n}}\bigl(\vg(s) -\vg(t_{n-1})\bigr)\, dW(s), \ve^{n-1}\biggr)\Bigr|,\\\nonumber
			D^n &:= \frac{2C_L^2}{\nu} \|\nab\vu^n\|^2_{\vL^2}\qquad\forall 1\leq n\leq M.
		\end{align*}
	\end{lemma}

	\begin{proof}
		Subtracting \eqref{eq_reforms} from \eqref{equu2.8a} we obtain
		\begin{align}
			\nonumber	\bigl(\ve^{n+1}-\ve^n,\vv\bigr) + \nu k\bigl(\nab\ve^{n+1},\nab\vv\bigr) &= \nu \int_{t_n}^{t_{n+1}}\bigl(\nab(\vu(s)-\vu(t_{n+1})),\nab\vv\bigr)\, ds \\
			&+ \int_{t_n}^{t_{n+1}}\bigl(\vu^{n+1}\cdot\nab\vu^{n+1} - \vu(s)\cdot\nab\vu(s),\vv\bigr)\, ds\\\nonumber
			&+\biggl(\int_{t_n}^{t_{n+1}}\bigl(\vg(s)-\vg(t_n)\bigr)\,dW(s),\vv\biggr).
		\end{align}
		Choosing $\vv = \ve^{n+1}$ yields
		\begin{align}\label{error_eq}
			\nonumber	\bigl(\ve^{n+1}-\ve^n,\ve^{n+1}\bigr) + \nu k\|\nab\ve^{n+1}\|^2_{\vL^2} &= \nu \int_{t_n}^{t_{n+1}}\bigl(\nab(\vu(s)-\vu(t_{n+1})),\nab\ve^{n+1}\bigr)\, ds \\\nonumber
			&+ \int_{t_n}^{t_{n+1}}\bigl(\vu^{n+1}\cdot\nab\vu^{n+1} - \vu(s)\cdot\nab\vu(s),\ve^{n+1}\bigr)\, ds\\
			&+\biggl(\int_{t_n}^{t_{n+1}}\bigl(\vg(s)-\vg(t_n)\bigr)\,dW(s),\ve^{n+1}\biggr)\\\nonumber
			&:= {\tt I + II + III}.
		\end{align}
		The left-hand side of \eqref{error_eq} can be easily controlled by using the formula $2(a,a-b) = \|a\|^2 - \|b\|^2 + \|a-b\|^2$. The right-hand side of \eqref{error_eq} can be bounded as follows.
		
		Using Cauchy-Schwarz inequality, we obtain
		\begin{align}\label{equu311}
			{\tt I} &\leq \nu \int_{t_n}^{t_{n+1}} \|\nab(\vu(t_{n+1}) - \vu(s))\|_{\vL^2}\|\nab\ve^{n+1}\|_{\vL^2}\, ds\\\nonumber
			&\leq \nu\int_{t_n}^{t_{n+1}}\|\nab(\vu(t_{n+1}) - \vu(s))\|^2_{\vL^2}\, ds + \frac{\nu k}{4}\|\nab\ve^{n+1}\|^2_{\vL^2}.
		\end{align}
		We rewrite the second term as
		\begin{align}\label{equu312}
			{\tt II} &= \int_{t_n}^{t_{n+1}}\bigl((\vu(s) - \vu(t_{n+1}))\cdot\nab \vu(s),\ve^{n+1}\bigr)\, ds \\\nonumber
			&\qquad+\int_{t_n}^{t_{n+1}} \bigl(\vu(t_{n+1})\cdot\nab(\vu(s) - \vu(t_{n+1})),\ve^{n+1}\bigr)\, ds \\\nonumber
			&\qquad+ \int_{t_n}^{t_{n+1}}\bigl(\ve^{n+1}\cdot\nab\vu^{n+1},\ve^{n+1}\bigr)\, ds \\\nonumber
			&\qquad+ \int_{t_n}^{t_{n+1}}\bigl(\vu^{n+1}\cdot\nab\ve^{n+1},\ve^{n+1}\bigr)\, ds\\\nonumber
			&=\int_{t_n}^{t_{n+1}}\bigl((\vu(s) - \vu(t_{n+1}))\cdot\nab \vu(s),\ve^{n+1}\bigr)\, ds \\\nonumber
			&\qquad+\int_{t_n}^{t_{n+1}} \bigl(\vu(t_{n+1})\cdot\nab(\vu(s) - \vu(t_{n+1})),\ve^{n+1}\bigr)\, ds \\\nonumber
			&\qquad+ \int_{t_n}^{t_{n+1}}\bigl(\ve^{n+1}\cdot\nab\vu^{n+1},\ve^{n+1}\bigr)\, ds\\\nonumber
			&= {\tt II_1 + II_2 +II_3}.
		\end{align}
		Here we have used the facts that $\bigl(\vu^{n+1}\cdot\nab\ve^{n+1},\ve^{n+1}\bigr) = 0$ and $\bigl(\ve^{n+1}\cdot\nab\vu^{n+1},\ve^{n+1}\bigr) = \bigl(\ve^{n+1}\cdot\nab\vu(t_{n+1}),\ve^{n+1}\bigr)$.
		In addition,
		\begin{align*}
			{\tt II_1} &\leq \int_{t_n}^{t_{n+1}} \|\vu(s) - \vu(t_{n+1})\|_{\vL^4}\|\vu(s)\|_{\vL^4}\|\nab\ve^{n+1}\|_{\vL^2}\, ds\\\nonumber
			&\leq \frac{2}{\nu}\int_{t_n}^{t_{n+1}} \|\vu(s) - \vu(t_{n+1})\|^2_{\vL^4}\|\vu(s)\|^2_{\vL^4} \, ds + \frac{\nu  k}{8}\|\nab\ve^{n+1}\|^2_{\vL^2}\\\nonumber
			&\leq \frac{2}{\nu}\int_{t_n}^{t_{n+1}} C_L^2\, \|\vu(s) - \vu(t_{n+1})\|^2_{\vH^1}\|\nab\vu(s)\|^2_{\vL^2} \, ds + \frac{\nu k}{8}\|\nab\ve^{n+1}\|^2_{\vL^2};
		\end{align*}
		\begin{align*}
			{\tt II_2} &\leq  \frac{2}{\nu}\int_{t_n}^{t_{n+1}} C_L^2\, \|\vu(s) - \vu(t_{n+1})\|^2_{\vH^1}\|\nab\vu(t_{n+1})\|^2_{\vL^2} \, ds + \frac{\nu k}{8}\|\nab\ve^{n+1}\|^2_{\vL^2};\\\nonumber
			{\tt II_3} &\leq k \|\ve^{n+1}\|^2_{\vL^4}\|\nab\vu(t_{n+1})\|_{\vL^2} \\\nonumber
			&\leq C_L k\|\nab\ve^{n+1}\|_{\vL^2}\|\ve^{n+1}\|_{\vL^2}\|\nab\vu(t_{n+1})\|_{\vL^2}\\\nonumber
			&\leq \frac{C_L^2}{\nu}\,k \, \|\ve^{n+1}\|^2_{\vL^2}\|\nab\vu^{n+1}\|^2_{\vL^2} + \frac{\nu k}{4}\|\nab\ve^{n+1}\|^2_{\vL^2}.
		\end{align*}
		Thus, \begin{align}\label{equu313}
			{\tt II} &\leq \frac{2}{\nu}\int_{t_n}^{t_{n+1}} C_L^2\, \|\vu(s) - \vu(t_{n+1})\|^2_{\vH^1}\|\nab\vu(t_{n+1})\|^2_{\vL^2} \, ds \\\nonumber
			&\qquad+ \frac{2}{\nu}\int_{t_n}^{t_{n+1}} C_L^2\, \|\vu(s) - \vu(t_{n+1})\|^2_{\vH^1}\|\nab\vu(t_{n+1})\|^2_{\vL^2} \, ds\\\nonumber
			&\qquad+ \frac{\nu k}{2}\|\nab\ve^{n+1}\|^2_{\vL^2} + \frac{C_L^2}{\nu}\, k\, \|\ve^{n+1}\|^2_{\vL^2}\|\nab\vu^{n+1}\|^2_{\vL^2}.
		\end{align}
		Now, we estimate the noise term as below
		\begin{align}\label{equu314}
			{\tt III} &= \biggl(\int_{t_n}^{t_{n+1}}\bigl(\vg(s) -\vg(t_n)\bigr)\, dW(s), \ve^{n+1} - \ve^n\biggr) \\\nonumber
			&\qquad+ \biggl(\int_{t_n}^{t_{n+1}}\bigl(\vg(s) -\vg(t_n)\bigr)\, dW(s), \ve^n\biggr)\\\nonumber
			&\leq \Bigl\|\int_{t_n}^{t_{n+1}}\bigl(\vg(s) - \vg(t_n)\bigr)\, dW(s)\Bigr\|^2_{\vL^2} + \frac14\|\ve^{n+1} - \ve^n\|^2_{\vL^2}\\\nonumber
			&\qquad+\biggl(\int_{t_n}^{t_{n+1}}\bigl(\vg(s) -\vg(t_n)\bigr)\, dW(s), \ve^n\biggr).
		\end{align}
		Substituting \eqref{equu311},\eqref{equu313} and \eqref{equu314} into \eqref{error_eq}, we obtain
		\begin{align}\label{equu315}
			\frac12 \bigl[\|\ve^{n+1}\|^2_{\vL^2} &- \|\ve^n\|^2_{\vL^2}\bigr] + \frac14 \|\ve^{n+1} - \ve^n\|^2_{\vL^2} + \frac{\nu k}{4}\|\nab\ve^{n+1}\|^2_{\vL^2} \\\nonumber
			&\leq \nu\int_{t_n}^{t_{n+1}}\|\nab(\vu(t_{n+1}) - \vu(s))\|^2_{\vL^2}\, ds \\\nonumber
			&\qquad+ \frac{2}{\nu}\int_{t_n}^{t_{n+1}} C_L^2\, \|\vu(s) - \vu(t_{n+1})\|^2_{\vH^1}\|\nab\vu(s)\|^2_{\vL^2} \, ds \\\nonumber
			&\qquad+ \frac{2}{\nu}\int_{t_n}^{t_{n+1}} C_L^2\, \|\vu(s) - \vu(t_{n+1})\|^2_{\vH^1}\|\nab\vu(t_{n+1})\|^2_{\vL^2} \, ds\\\nonumber
			&\qquad+\frac{C_L^2}{\nu}\, k\, \|\ve^{n+1}\|^2_{\vL^2}\|\nab\vu^{n+1}\|^2_{\vL^2}\\\nonumber
			&\qquad+ \Bigl\|\int_{t_n}^{t_{n+1}}\bigl(\vg(s) - \vg(t_n)\bigr)\, dW(s)\Bigr\|_{\vL^2}^2 \\\nonumber
			&\qquad+ \biggl(\int_{t_n}^{t_{n+1}}\bigl(\vg(s) -\vg(t_n)\bigr)\, dW(s), \ve^n\biggr).
		\end{align}
		Next, lowering the index $n$ in \eqref{equu315} by $1$ and applying the summation operator $ \sum_{n=1}^{\ell}$ to both sides of \eqref{equu315} for any $1\leq \ell \leq M$, we obtain
		\begin{align}\label{equu316}
			&\|\ve^{\ell}\|^2_{\vL^2} + \frac{1}{2}\sum_{n=1}^{\ell} \|\ve^{n} - \ve^{n-1}\|^2_{\vL^2} + \frac{\nu k}{2}\sum_{n=1}^{\ell} \|\nab\ve^{n}\|^2_{\vL^2} \\\nonumber
			&\leq 2\nu\sum_{n=1}^{\ell}\int_{t_{n-1}}^{t_{n}}\|\nab(\vu(t_{n}) - \vu(s))\|^2_{\vL^2}\, ds \\\nonumber
			&\qquad+ \frac{4}{\nu}\sum_{n=1}^{\ell}\int_{t_{n-1}}^{t_{n}} C_L^2\, \|\vu(s) - \vu(t_{n})\|^2_{\vH^1}\|\nab\vu(s)\|^2_{\vL^2} \, ds \\\nonumber
			&\qquad+ \frac{4}{\nu}\sum_{n=1}^{\ell}\int_{t_{n-1}}^{t_{n}} C_L^2\, \|\vu(s) - \vu(t_{n})\|^2_{\vH^1}\|\nab\vu(t_{n})\|^2_{\vL^2} \, ds\\\nonumber
			&\qquad+\frac{2C_L^2}{\nu}\, k\,\sum_{n=1}^{\ell} \|\ve^{n}\|^2_{\vL^2}\|\nab\vu^{n}\|^2_{\vL^2}\\\nonumber
			&\qquad+ 2\sum_{n=1}^{\ell}\Bigl\|\int_{t_{n-1}}^{t_{n}}\bigl(\vg(s) - \vg(t_{n-1})\bigr)\, dW(s)\Bigr\|_{\vL^2}^2\\\nonumber
			&\qquad+ 2\sum_{n=1}^{\ell}\biggl(\int_{t_{n-1}}^{t_{n}}\bigl(\vg(s) -\vg(t_{n-1})\bigr)\, dW(s), \ve^{n-1}\biggr)\\\nonumber
			&\leq 2\nu\sum_{n=1}^{M}\int_{t_{n-1}}^{t_{n}}\|\nab(\vu(t_{n}) - \vu(s))\|^2_{\vL^2}\, ds \\\nonumber
			&\qquad+ \frac{4}{\nu}\sum_{n=1}^{M}\int_{t_{n-1}}^{t_{n}} C_L^2\, \|\vu(s) - \vu(t_{n})\|^2_{\vH^1}\|\nab\vu(s)\|^2_{\vL^2} \, ds \\\nonumber
			&\qquad+ \frac{4}{\nu}\sum_{n=1}^{M}\int_{t_{n-1}}^{t_{n}} C_L^2\, \|\vu(s) - \vu(t_{n})\|^2_{\vH^1}\|\nab\vu(t_{n})\|^2_{\vL^2} \, ds\\\nonumber
			&\qquad+\frac{2C_L^2}{\nu}\, k\,\sum_{n=1}^{\ell-1} \|\ve^{n}\|^2_{\vL^2}\|\nab\vu^{n}\|^2_{\vL^2} + \frac{2C_L^2k}{\nu}\max_{1\leq \ell \leq M}\|\ve^{\ell}\|^2_{\vL^2}\|\nab\vu^{\ell}\|^2_{\vL^2}\\\nonumber
			&\qquad+ 2\sum_{n=1}^{M}\Bigl\|\int_{t_{n-1}}^{t_{n}}\bigl(\vg(s) - \vg(t_{n-1})\bigr)\, dW(s)\Bigr\|_{\vL^2}^2\\\nonumber
			&\qquad+ 2\max_{1\leq \ell \leq M}\Bigl|\sum_{n=1}^{\ell}\biggl(\int_{t_{n-1}}^{t_{n}}\bigl(\vg(s) -\vg(t_{n-1})\bigr)\, dW(s), \ve^{n-1}\biggr)\Bigr|\\\nonumber
			&= A + B + k\sum_{n = 1}^{\ell-1} D^n\|\ve^n\|^2_{\vL^2},
		\end{align}
		where $A,B$ are constants in $\ell$ and are defined by
		\begin{align}\label{equu317}
			A &:= 2\nu\sum_{n=1}^{M}\int_{t_{n-1}}^{t_{n}}\|\nab(\vu(t_{n}) - \vu(s))\|^2_{\vL^2}\, ds \\\nonumber
			&\qquad+ \frac{4}{\nu}\sum_{n=1}^{M}\int_{t_{n-1}}^{t_{n}} C_L^2\, \|\vu(s) - \vu(t_{n})\|^2_{\vH^1}\|\nab\vu(s)\|^2_{\vL^2} \, ds \\\nonumber
			&\qquad+ \frac{4}{\nu}\sum_{n=1}^{M}\int_{t_{n-1}}^{t_{n}} C_L^2\, \|\vu(s) - \vu(t_{n})\|^2_{\vH^1}\|\nab\vu(t_{n})\|^2_{\vL^2} \, ds\\\nonumber
			&\qquad+ \frac{2C_L^2}{\nu} k \max_{1\leq \ell \leq M}\|\ve^{\ell}\|^2_{\vL^2}\|\nab\vu^{\ell}\|^2_{\vL^2}, \\\nonumber
			B &:= 2\sum_{n=1}^{M}\Bigl\|\int_{t_{n-1}}^{t_{n}}\bigl(\vg(s) - \vg(t_{n-1})\bigr)\, dW(s)\Bigr\|_{\vL^2}^2\\\nonumber
			&\qquad+ 2\max_{1\leq \ell \leq M}\Bigl|\sum_{n=1}^{\ell}\biggl(\int_{t_{n-1}}^{t_{n}}\bigl(\vg(s) -\vg(t_{n-1})\bigr)\, dW(s), \ve^{n-1}\biggr)\Bigr|,\\\nonumber
			D^n &:= \frac{2C_L^2}{\nu} \|\nab\vu^n\|^2_{\vL^2}\qquad\forall 1\leq n\leq M.
		\end{align}
		Finally, applying the discrete Gronwall inequality	to \eqref{equu316}, we obtain
		\begin{align}\label{equu318}
			\|\ve^{\ell}\|^2_{\vL^2} &+  \sum_{n=1}^{M} \|\ve^{n} - \ve^{n-1}\|^2_{\vL^2} + {\nu k}\sum_{n=1}^{M} \|\nab\ve^{n}\|^2_{\vL^2} \\\nonumber
			&\leq \bigl(A + B\bigr)\exp\biggl(k\sum_{n=1}^{\ell-1} D^n\biggr) \qquad\mP-\mbox{a.s.}
		\end{align}
		The proof is complete by taking $\displaystyle\max_{1\leq \ell \leq M}$ on both sides.
	\end{proof}
	
	Next, we estimate $ \mE\biggl[\bigl(A + B\bigr)^{\frac{q}{2}}\exp\biggl(\frac{kq}{2}\sum_{n=1}^{M} D^n\biggr)\biggr]$ in the following lemma to extract the convergent rates as stated in Theorem \ref{theorem_semi_chapter5}. 
	
	\smallskip
	\begin{lemma}\label{lemma3.5}
		Assume that (G2) satisfies with $\sigma = \frac{2C_L^2}{\nu}$. Then, under the assumptions of Theorem \ref{theorem_semi_chapter5} and Lemma \ref{lemma3.4}, there holds
		\begin{align}
			\biggl(\mE\biggl[\bigl(A + B\bigr)^{\frac{q}{2}}\exp\biggl(\frac{kq}{2}\sum_{n=1}^{M} D^n\biggr)\biggr]\biggr)^{\frac1q} \leq C_q\, k^{\frac12-\gamma}.
		\end{align}
		Here, $A,B,D^n$ are defined as in \eqref{equu317} and $0<\gamma<\frac12$. $C_q$ is the same constant as in Theorem \ref{theorem_semi_chapter5} which depends on $\nu, q, \vg, D_T, \vu_0$.
	\end{lemma}
	
	\begin{proof} 
		By using the simple inequality $(a+b)^p \leq 2^{p-1}(a^p + b^p)$ for $1\leq p <\infty$ we have
		\begin{align}\label{equu3.40}
			\mE\biggl[\bigl(A + B\bigr)^{\frac{q}{2}}\exp\biggl(\frac{kq}{2}\sum_{n=1}^{M} D^n\biggr)\biggr] &\leq C_q\mE\biggl[A^{\frac{q}{2}}\exp\Bigl(\frac{kq}{2}\sum_{n = 1}^M D^n\Bigr)\biggr] \\\nonumber
			&\qquad+ C_q\mE\biggl[B^{\frac{q}{2}}\exp\Bigl(\frac{kq}{2}\sum_{n = 1}^M D^n\Bigr)\biggr]\\\nonumber
			&\leq C_q\bigl(\mE[A^q]\bigr)^{\frac12}\Bigl(\mE\Bigl[\exp\Bigl(qk\sum_{n = 1}^{M}D^n\Bigr)\Bigr]\Bigr)^{\frac12}\\\nonumber
			&\quad+C_q\bigl(\mE[B^q]\bigr)^{\frac12}\Bigl(\mE\Bigl[\exp\Bigl(qk\sum_{n = 1}^{M}D^n\Bigr)\Bigr]\Bigr)^{\frac12}.
		\end{align}
		Therefore, there are three factors in \eqref{equu3.40} which we need to estimate, namely, $ \mE[A^q],\, \mE[B^q]$ and $ \mE\biggl[\exp\biggl(qk\sum_{n=1}^{M}D^n\biggr)\biggr]$.
		
		\begin{enumerate}[(a)]
			\item To estimate $\mE[A^q]$, by \eqref{equu317} and the discrete Cauchy-Schwarz inequality, we have
			\begin{align*}
				\mE[A^q] &= \mE\biggl[\biggl(2\nu\sum_{n=1}^{M}\int_{t_{n-1}}^{t_{n}}\|\nab(\vu(t_{n}) - \vu(s))\|^2_{\vL^2}\, ds \\\nonumber
				&\qquad+ \frac{4}{\nu}\sum_{n=1}^{M}\int_{t_{n-1}}^{t_{n}} C_L^2\, \|\vu(s) - \vu(t_{n})\|^2_{\vH^1}\|\nab\vu(s)\|^2_{\vL^2} \, ds \\\nonumber
				&\qquad+ \frac{4}{\nu}\sum_{n=1}^{M}\int_{t_{n-1}}^{t_{n}} C_L^2\, \|\vu(s) - \vu(t_{n})\|^2_{\vH^1}\|\nab\vu(t_{n})\|^2_{\vL^2} \, ds\\\nonumber
				&\qquad+ \frac{2C_L^2}{\nu} k \max_{1\leq \ell \leq M}\|\ve^{\ell}\|^2_{\vL^2}\|\nab\vu^{\ell}\|^2_{\vL^2}	\biggr)^q\biggr]\\\nonumber
				&\leq C_q\mE\biggl[\biggl(\sum_{n=1}^{M}\int_{t_{n-1}}^{t_{n}}\|\nab(\vu(t_{n}) - \vu(s))\|^2_{\vL^2}\, ds\biggr)^q\biggr]\\\nonumber
				&\qquad+ C_q\mE\biggl[\biggl(\sup_{s\in [0,T]}\|\nab\vu(s)\|^2_{\vL^2}\sum_{n=1}^{M}\int_{t_{n-1}}^{t_{n}} \, \|\vu(s) - \vu(t_{n})\|^2_{\vH^1} \, ds\biggr)^q\biggr]\\\nonumber
				&\qquad+ C_qk^q\mE\bigl[\max_{1\leq n \leq M}\|\ve^n\|^q_{\vL^2}\|\nab\vu^n\|^q_{\vL^2}\bigr]\\\nonumber
				& = A_1 +A_2 +A_3.
			\end{align*}
			
			By using Cauchy-Schwarz inequality and \eqref{equu2.20a}, we have
			\begin{align*}
				A_1 &= C_q\mE\biggl[\biggl(\sum_{n=1}^{M}\int_{t_{n-1}}^{t_{n}}\|\nab(\vu(t_{n}) - \vu(s))\|^{2}_{\vL^2}\, ds\biggr)^q\biggr]\\\nonumber
				&\leq C_q\mE\biggl[\sum_{n=1}^{M}\int_{t_{n-1}}^{t_{n}}\|\nab(\vu(t_{n}) - \vu(s))\|^{2q}_{\vL^2}\, ds\biggr]\\\nonumber
				&\leq C_qk^{2q(\frac12-\gamma)}.
			\end{align*}
			By using Cauchy-Schwarz inequality, \eqref{stability_means} and \eqref{equu2.20a}, we also obtain
			\begin{align*}
				A_2 &= C_q\mE\biggl[\biggl(\sup_{s\in [0,T]}\|\nab\vu(s)\|^2_{\vL^2}\sum_{n=1}^{M}\int_{t_{n-1}}^{t_{n}} \, \|\vu(s) - \vu(t_{n})\|^2_{\vH^1} \, ds\biggr)^q\biggr]\\\nonumber
				&\leq C_q\mE\biggl[\sup_{s\in [0,T]}\|\nab\vu(s)\|^{2q}_{\vL^2}\sum_{n=1}^M\int_{t_{n-1}}^{t_n}\|\vu(s) - \vu(t_n)|^{2q}_{\vH^1}\, ds\biggr]\\\nonumber
				&\leq C_q\bigl(\mE\bigl[\sup_{s\in [0,T]}\|\nab\vu(s)\|^{4q}_{\vL^2}\bigr]\bigr)^{\frac12}\biggl(\mE\biggl[\sum_{n=1}^M\int_{t_{n-1}}^{t_n}\|\vu(s) - \vu(t_n)|^{4q}_{\vH^1}\, ds\biggr]\biggr)^{\frac12}\\\nonumber
				&\leq C_qk^{2q(\frac12-\gamma)}.
			\end{align*}
			In addition, by using \eqref{stability_means}, $A_3 \leq C_q k^q$. In summary, $\mE\bigl[A^q\bigr] \leq C_q k^{2q(\frac12-\gamma)}$.
			
			\smallskip
			\item To estimate $\mE[B^q]$, by \eqref{equu317}, we have
			\begin{align*}
				\mE[B^q] &= \mE\biggl[\biggl(2\sum_{n=1}^{M}\Bigl\|\int_{t_{n-1}}^{t_{n}}\bigl(\vg(s) - \vg(t_{n-1})\bigr)\, dW(s)\Bigr\|_{\vL^2}^2\\\nonumber
				&\qquad+ 2\max_{1\leq \ell \leq M}\Bigl|\sum_{n=1}^{\ell}\Bigl(\int_{t_{n-1}}^{t_{n}}\bigl(\vg(s) -\vg(t_{n-1})\bigr)\, dW(s), \ve^{n-1}\Bigr)\Bigr|\biggr)^q\biggr]\\\nonumber
				&\leq C_q\mE\biggl[\biggl(\sum_{n = 1}^M\Bigl\|\int_{t_{n-1}}^{t_{n}}\bigl(\vg(s) - \vg(t_{n-1})\bigr)\, dW(s)\Bigr\|^2_{\vL^2}\biggr)^q\biggr]\\\nonumber
				&\qquad + C_q\mE\biggl[\max_{1\leq \ell \leq M}\Bigl|\sum_{n=1}^{\ell}\Bigl(\int_{t_{n-1}}^{t_{n}}\bigl(\vg(s) -\vg(t_{n-1})\bigr)\, dW(s), \ve^{n-1}\Bigr)\Bigr|^q\biggr]\\\nonumber
				&= B_1 + B_2.
			\end{align*}
			Next, by using the discrete H\"older inequality, \eqref{ito4} and (G1) we have
			
			\begin{align*}
				B_1 &= C_q\mE\biggl[\biggl(\sum_{n = 1}^M\Bigl\|\int_{t_{n-1}}^{t_{n}}\bigl(\vg(s) - \vg(t_{n-1})\bigr)\, dW(s)\Bigr\|^2_{\vL^2}\biggr)^q\biggr]\\\nonumber
				&\leq C_q M^{q-1}\mE\biggl[\sum_{n = 1}^M\Bigl\|\int_{t_{n-1}}^{t_n}\bigl(\vg(s) - \vg(t_{n-1})\bigr)\, dW(s)\Bigr\|^{2q}_{\vL^2}\biggr]\\\nonumber
				&\leq C_qM^{q-1}\mE\biggl[\sum_{n = 1}^M\int_{t_{n-1}}^{t_n} \|\vg(s) - \vg(t_{n-1})\|^{2q}_{\vL^2}\, ds\biggr]\\\nonumber&\leq C_q k^{q+1}.
			\end{align*}
			
			To estimate $B_2$, we use the Burkholder--Davis--Gundy inequality and (G1) to get
			\begin{align*}
				B_2 &= C_q\mE\biggl[\max_{1\leq \ell \leq M}\Bigl|\sum_{n=1}^{\ell}\Bigl(\int_{t_{n-1}}^{t_{n}}\bigl(\vg(s) -\vg(t_{n-1})\bigr)\, dW(s), \ve^{n-1}\Bigr)\Bigr|^q\biggr]\\\nonumber
				&\leq C_q\mE\biggl[\Bigl(\sum_{n = 1}^M\int_{t_{n-1}}^{t_n}\|\vg(s) - \vg(t_{n-1})\|^2_{\vL^2}\|\ve^{n-1}\|^{2}_{\vL^2}\, ds\Bigr)^{q/2}\biggr]\\\nonumber
				&\leq C_qk^q \mE\biggl[\Bigl(k\sum_{n = 1}^M\|\ve^{n-1}\|^2_{\vL^2}\Bigr)^{q/2}\biggr]\leq C_qk^q.
			\end{align*}
			In summary, $$\mE\bigl[B^q\bigr] \leq C_q k^q. $$

			\smallskip
			\item To estimate $ \mE\biggl[\exp\biggl(qk\sum_{n=1}^{M}D^n\biggr)\biggr]$ by \eqref{equu317} and the discrete Jensen inequality, we have
			\begin{align}\label{equu3.23}
			\nonumber 	\mE\biggl[\exp\biggl(qk\sum_{n=1}^{M}D^n\biggr)\biggr] &\leq \mE\biggl[\exp\biggl(\frac{2qC_L^2}{\nu}k\sum_{n=0}^{M}\|\nab\vu^n\|^2_{\vL^2}\biggr)\biggr]\\\nonumber
				&\leq \mE\biggl[k\sum_{n=0}^M \exp\biggl(\frac{2qC_L^2}{\nu}\|\nab\vu^n\|^2_{\vL^2}\biggr)\biggr]\\
				&\leq k \sum_{n=0}^M \mE\biggl[\max_{0\leq n \leq M}\exp\biggl(\frac{2qC_L^2}{\nu}\|\nab\vu^n\|^2_{\vL^2}\biggr)\biggr]
				 \leq CT. 
			\end{align}
			Here, we have used Lemma \ref{lemma_exp_discrete} with $\sigma = \frac{2qC^2_L}{\nu}$ to get the last inequality of \eqref{equu3.23} .

		\end{enumerate}
		
		Finally, combining (a), (b), (c) and \eqref{equu3.40}, we get
		\begin{align}\label{equu325}
			\mE\biggl[\bigl(A + B\bigr)^{\frac{q}{2}}\exp\biggl(\frac{kq}{2}\sum_{n=1}^{M} D^n\biggr)\biggr]\leq C_q k^{(\frac12-\gamma)q}.
		\end{align}
		The proof is complete.
	\end{proof}
	
	We are ready to state our first pathwise error estimate for the velocity approximation,
	such an estimate has not been obtained before in the literature. 
	
	\begin{theorem}\label{thm_pathwise_semi}
		Assume that the assumptions of Theorem \ref{theorem_semi_chapter5} hold. Let $2 < q < \infty$ and $0<\gamma < \frac{1}{2}$ such that $\frac12 - \gamma - \frac{1}{q} >0$. Then, for $0 < \gamma_1 < \frac12 -  \gamma - \frac{1}{q}$, there exists a random variable $K = K(\omega)$ with $\mE\bigl[|K|^{q}\bigr] <\infty$ such that there holds $\mP$-a.s.
		\begin{align}\label{equu3.41}
			\max_{1\leq n \leq M}\|\vu(t_n) - \vu^n\|_{\vL^2}  + \Bigl(\nu k\sum_{n=1}^M\|\nab\bigl(\vu(t_n) - \vu^n\bigr)\|^2_{\vL^2}\Bigr)^{\frac12} \leq K\, k^{\gamma_1}.
		\end{align}
	\end{theorem}

	\begin{proof}
		\eqref{equu3.41} is an immediate consequence of Theorem \ref{theorem_semi_chapter5} and Kolmogorov Theorem \ref{kolmogorov}.
	\end{proof}
	
	\subsection{High moment and pathwise error estimates for the pressure approximation} 
	We state and prove the following pathwise error estimate for the pressure approximation. Such a result has not been obtained before in the literature. 
	\begin{theorem}\label{semi_pressure_error}
		Assume that the assumptions of Theorem \ref{theorem_semi_chapter5} hold. Then we have
		\begin{align}
			\biggl(\mE\biggl[\Bigl\|P(t_m) - k\sum_{n=1}^m p^n\Bigr\|^q_{\vL^2}\biggr]\biggr)^{\frac1q} \leq C_qk^{\frac12 - \gamma}.
		\end{align}
	\end{theorem}
	\begin{proof} The proof is based on the well-known inf-sup (LBB) condition which is quoted below. There exists a constant $\beta >0$ such that
		\begin{align}
			\beta \|\xi\|_{L^2} \leq \sup_{\pphi \in \vH^1_{per}(D)}\frac{\bigl(\xi,\div \pphi\bigr)}{\| \nab\pphi \|_{\vL^2}}\,\,\qquad\qquad\forall \xi \in L_{per}^2(D).
		\end{align}
		Set
		\begin{align*}
			\ve^m := \vu(t_m) - \vu^m \mbox{ and }\E^m := P(t_m) - P^m, \mbox{ where } P^m := k\sum_{n = 1}^m p^{n}.
		\end{align*} 
		Summing \eqref{equu3.1} from $1$ to $m$ we obtain
		\begin{align}\label{eq_summation_semi}
			\bigl(\vu^m,\pphi\bigr) &+ \nu k\sum_{n = 1}^m\bigl(\nab\vu^n,\nab\pphi\bigr) + k\sum_{n = 1}^m\bigl(\vu^n\cdot\nab\vu^n,\pphi\bigr) - \bigl(P^m,\div \pphi\bigr) \\\nonumber
			&=\bigl(\vu^0,\pphi\bigr) + \sum_{n = 1}^m\bigl(\vg^{n-1}\Delta W_n,\pphi\bigr)\,\,\qquad\qquad\forall\pphi\in \vH^1_{per}(D).
		\end{align}
		Subtracting \eqref{eq_summation_semi} from \eqref{equu2.10a}, we obtain the following error equations:
		\begin{align}
			\bigl(\ve^m,\pphi\bigr) &+ \Bigl(\nu k\sum_{n=1}^m \nab\ve^n,\nab\pphi\Bigr) + \nu k \sum_{n = 1}^m\int_{t_{n-1}}^{t_n}\bigl(\nab(\vu(s)-\vu(t_n)),\nab\pphi\bigr)\, ds \\\nonumber
			&+ k\sum_{n = 1}^m \bigl(\vu(t_n)\cdot\nab\vu(t_n) - \vu^n\cdot\nab\vu^n, \pphi\bigr) \\\nonumber
			&+ \sum_{n=1}^m \int_{t_{n-1}}^{t_n}\bigl(\vu(s)\cdot\nab\vu(s) - \vu(t_n)\cdot\nab\vu(t_n),\pphi\bigr)\, ds - \bigl(\E^m,\div \pphi\bigr)\\\nonumber
			&= \Bigl(\sum_{n = 1}^m\int_{t_{n-1}}^{t_n}(\vg(s) - \vg^{n-1})\, dW(s),\pphi\Bigr).
		\end{align}
		Therefore,
		\begin{align}\label{equu3.50}
			\nonumber	\bigl(\E^m,\div \pphi\bigr) &= \bigl(\ve^m,\pphi\bigr) + \Bigl(\nu k\sum_{n=1}^m \nab\ve^n,\nab\pphi\Bigr) + \nu  \sum_{n = 1}^m\int_{t_{n-1}}^{t_n}\bigl(\nab(\vu(s)-\vu(t_n)),\nab\pphi\bigr)\, ds \\\nonumber
			&\qquad+  k\sum_{n = 1}^m \bigl(\vu(t_n)\cdot\nab\vu(t_n) - \vu^n\cdot\nab\vu^n, \pphi\bigr) \\
			&\qquad+ \sum_{n=1}^m \int_{t_{n-1}}^{t_n}\bigl(\vu(s)\cdot\nab\vu(s) - \vu(t_n)\cdot\nab\vu(t_n),\pphi\bigr)\, ds\\\nonumber
			&\qquad+ \Bigl(\sum_{n = 1}^m\int_{t_{n-1}}^{t_n}(\vg(t_{n-1}) - \vg(s))\, dW(s),\pphi\Bigr)\\\nonumber
			&:= {\tt I + II + III + IV + V +VI}.
		\end{align}
		Next, we use Cauchy-Schwarz and Poincar\'e inequalities to estimate each term on the right-hand side of \eqref{equu3.50}.
		\begin{align*} 
			{\tt I + II} &\leq C\|\ve^m\|_{\vL^2}\|\nab\pphi\|_{\vL^2} + \nu k \sum_{n = 1}^m\|\nab\ve^n\|_{\vL^2}\|\nab\pphi\|_{\vL^2}. \\\nonumber
			{\tt III} &\leq \nu \sum_{n=1}^m\int_{t_{n-1}}^{t_n}\|\nab(\vu(s) - \vu(t_n))\|_{\vL^2}\, ds \|\nab\pphi\|_{\vL^2}.\\\nonumber
			{\tt IV} &= k\sum_{n = 1}^m \bigl(\vu(t_n)\cdot\nab\vu(t_n) - \vu^n\cdot\nab\vu^n, \pphi\bigr)\\\nonumber
			&=k\sum_{n = 1}^m \bigl(\vu(t_n)\cdot\nab\vu(t_n) - \vu^n\cdot\nab\vu(t_n), \pphi\bigr) \\\nonumber
			&\qquad+ k\sum_{n = 1}^m \bigl(\vu^n\cdot\nab\vu(t_n) - \vu^n\cdot\nab\vu^n, \pphi\bigr) \\\nonumber
			&= k\sum_{n = 1}^m \bigl(\ve^n\cdot\nab\vu(t_n), \pphi\bigr) + k\sum_{n = 1}^m \bigl(\vu^n\cdot\nab\ve^n, \pphi\bigr)\\\nonumber
			&\leq k\sum_{n = 1}^m \|\ve^n\|_{\vL^4}\|\vu(t_n)\|_{\vL^4}\|\nab\pphi\|_{\vL^2} + k\sum_{n = 1}^m\|\vu^n\|_{\vL^4}\|\ve^n\|_{\vL^4}\|\nab\pphi\|_{\vL^2}. \\
			{\tt V} &= \sum_{n=1}^m \int_{t_{n-1}}^{t_n}\bigl(\vu(s)\cdot\nab\vu(s) - \vu(t_n)\cdot\nab\vu(s),\pphi\bigr)\, ds \\\nonumber
			&\qquad+ \sum_{n=1}^m \int_{t_{n-1}}^{t_n}\bigl(\vu(t_n)\cdot\nab\vu(s) - \vu(t_n)\cdot\nab\vu(t_n),\pphi\bigr)\, ds\\\nonumber
			&\leq \sum_{n = 1}^m\int_{t_{n-1}}^{t_n}\|\vu(s) - \vu(t_n)\|_{\vL^4}\|\vu(s)\|_{\vL^4}\|\nab\pphi\|_{\vL^2}\, ds \\\nonumber
			&\qquad+ \sum_{n = 1}^m\int_{t_{n-1}}^{t_n}\|\vu(t_n)\|_{\vL^4}\|\vu(s) - \vu(t_n)\|_{\vL^4}\|\nab\pphi\|_{\vL^2}\, ds.\\\nonumber
			{\tt VI} &\leq C\Bigl\|\sum_{n = 1}^m\int_{t_{n-1}}^{t_n}(\vg(s) -\vg(t_{n-1}))\, dW(s)\Bigr\|_{\vL^2}\|\nab\pphi\|_{\vL^2}.
		\end{align*}
		Substituting the above estimates into \eqref{equu3.50} and using the inf-sup condition yield
		\begin{align}\label{equu3.57}
			\nonumber	\beta\|\E^m\|_{L^2} &\leq \sup_{\pphi\in \vH^1_{per}(D)}\frac{\bigl(\E^m,\div \pphi\bigr)}{\|\nab\pphi\|_{\vL^2}}\\\nonumber
			&\leq C\|\ve^m\|_{\vL^2} + \nu k \sum_{n = 1}^m\|\nab\ve^n\|_{\vL^2} +\nu \sum_{n=1}^m\int_{t_{n-1}}^{t_n}\|\nab(\vu(s) - \vu(t_n))\|_{\vL^2}\, ds \\\nonumber
			&\qquad + k\sum_{n = 1}^m \|\ve^n\|_{\vL^4}(\|\vu(t_n)\|_{\vL^4} + \|\vu^n\|_{\vL^4})\\
			&\qquad + \sum_{n = 1}^m\int_{t_{n-1}}^{t_n}\|\vu(s) - \vu(t_n)\|_{\vL^4}\bigl(\|\vu(s)\|_{\vL^4} + \|\vu(t_n)\|_{\vL^4}\bigr)\, ds\\\nonumber
			&\qquad + C\Bigl\|\sum_{n = 1}^m\int_{t_{n-1}}^{t_n}(\vg(s) -\vg(t_{n-1}))\, dW(s)\Bigr\|_{\vL^2}.
		\end{align}
		Taking the expectation on \eqref{equu3.57} and then using Theorem \ref{theorem_semi_chapter5}, Lemmas \ref{lemma2.2}--\ref{stability_means} and It\^o isometry, we obtain
		\begin{align*}
			\nonumber	\beta^q\mE\bigl[\|\E^m\|^q_{L^2}\bigr] &\leq C_q\mE[\|\ve^m\|^q_{\vL^2}] + C_q\mE\Bigl[\Bigl(\nu k\sum_{n=1}^m\|\nab\ve^n\|^2_{\vL^2}\Bigr)^{\frac{q}{2}}\Bigr] \\\nonumber
			&+ C_q \mE\Bigl[\sum_{n = 1}^m\int_{t_{n-1}}^{t_n}\|\nab(\vu(s) - \vu(t_n))\|^q_{\vL^2}\Bigr]\\\nonumber
			&+ C_q\mE\Bigl[\Bigl(k\sum_{n = 1}^m \|\nab\ve^n\|_{\vL^2}\bigl(\|\nab\vu(t_n)\|_{\vL^2} + \|\nab\vu^n\|_{\vL^2}\bigr)\Bigr)^{q}\Bigr]\\\nonumber
			&+ C_q\mE\Bigl[\sum_{n = 1}^m\int_{t_{n-1}}^{t_n}\|\nab(\vu(s) - \vu(t_n))\|^q_{\vL^2}\bigl(\|\nab\vu(s)\|^q_{\vL^2} + \|\nab\vu(t_n)\|^q_{\vL^2}\bigr)\, ds\Bigr]\\
			&+ C_q\mE\Bigl[\Bigl(\sum_{n = 1}^M\int_{t_{n-1}}^{t_n}\|\vg(s) - \vg(t_{n-1})\|^2_{\vL^2}\, ds\Bigr)^{\frac{q}{2}}\Bigr]\\\nonumber
			&\leq C_q k^{(\frac12 - \gamma)q} \\\nonumber
			&+ C_q\mE\Bigl[\Bigl(k\sum_{n = 1}^m\|\nab\ve^n\|^2_{\vL^2}\Bigr)^{\frac{q}{2}}\Bigl(k\sum_{n = 1}^m\bigl(\|\nab\vu(t_n)\|_{\vL^2} + \|\nab\vu^n\|_{\vL^2}\bigr)^2\Bigr)^{\frac{q}{2}}\Bigr]\\\nonumber
			&+ C_q\bigl(\mE\bigl[\max_{s \in [0,T]}\|\nab\vu(s)\|^{2q}_{\vL^2}\bigr]\bigr)^{\frac12}\sum_{n = 1}^m\int_{t_{n-1}}^{t_n}\bigl(\mE\bigl[\|\nab(\vu(s) - \vu(t_n))\|^{2q}_{\vL^2}\, \bigr]\bigr)^{1/2}ds\\\nonumber
			&\leq C_q k^{(\frac12 - \gamma)q} \\\nonumber
			&+ C_q\Bigl(\mE\Bigl[\Bigl(k\sum_{n = 1}^m\|\nab\ve^n\|^2_{\vL^2}\Bigr)^{q}\Bigr]\Bigr)^{\frac12}\Bigl(\mE\Bigl[\Bigl(k\sum_{n = 1}^m\bigl(\|\nab\vu(t_n)\|_{\vL^2} + \|\nab\vu^n\|_{\vL^2}\bigr)^2\Bigr)^{q}\Bigr]\Bigr)^{\frac12}\\\nonumber
			&\leq C_q k^{(\frac12 - \gamma)q}.
		\end{align*}
		Finally, the proof is complete after  dividing  both sides by $\beta^q$.
	\end{proof}

	Next, we state a pathwise error estimate for the pressure approximation. To the best of our knowledge, this is the first pathwise convergence result for  the pressure approximation.
	\begin{theorem}\label{thm_semi_pathwise_pressure}
		Assume the assumptions of Theorem \ref{semi_pressure_error} hold. Let $2 < q < \infty$ and $0 < \gamma < \frac12$ such that $\frac12 - \gamma - \frac{1}{q} > 0$. Then, for $0 < \gamma_1 < \frac12 - \gamma - \frac{1}{q}$, there exists a random variable $K = K(\omega)$ with $\mE\bigl[|K|^{q}\bigr] <\infty$ such that for all $1\leq \ell \leq M$, there holds $\mP$-a.s.
		\begin{align}\label{equu3.52}
			\biggl\| P(t_{\ell}) - k\sum_{n=1}^{\ell} p^n\biggr\|_{L^2} \leq K\, k^{\gamma_1}.
		\end{align}
	\end{theorem}
	\begin{proof}
		The assertion follows immediately from an application of Theorem \ref{kolmogorov} using the high moment error estimate of Theorem \ref{semi_pressure_error}.
	\end{proof}
	
	\section{Mixed finite element discretization in space}\label{section5.4}
	 In this section, we consider the spatial approximation for the Algorithm 1 by using finite element method. Let $\mathcal{T}_h$ be a quasi-uniform triangular mesh of the given domain $D \subset \mathbb{R}^2$ with mesh size $h > 0$. We introduce the following Taylor-Hood mixed finite element spaces:
	\begin{align*}
		\mH_h &:= \bigl\{\vv_h \in C(\overline{D}) \cap \vH^1_{per}(D);\,\vv_h \in  [\mathcal{P}_2(K)]^2\, \,\forall K \in \mathcal{T}_h\bigr\},\\
		L_h &:= \bigl\{q_h \in L^2_{per}(D)/\mathbb{R};\, q_h \in \mathcal{P}_1(K)\,\,\forall K \in \mathcal{T}_h\bigr\},
	\end{align*}  
	where $\mathcal{P}_i(K)$ denotes the space of all polynomials on $K$ of degree at most $i$. It's well-known \cite{Brezzi_Fortin91} that the above pair $(\mH_h,L_h)$   satisfies the Ladyzh\v{e}nskaja-Babu\v{s}ka-Brezzi (LBB) condition,  
	namely
	\begin{align}
		\sup_{\pphi_h \in \mH_h} \frac{\bigl(\div \pphi_h,q_h\bigr)}{\|\nab\pphi_h\|_{\vL^2}} \geq C\|q_h\|_{L^2}\qquad\forall q_h\in L_h,
	\end{align}
	where the constant $C>0$ is independent of $h$ and $k$. 
	
	Next, we define the space of weakly divergent-free velocity field as follows:
	\begin{align*}
		\mV_h := \bigl\{\pphi_h \in \mH_h;\, \bigl(\div \pphi_h, q_h\bigr) =0\qquad\forall q_h \in L_h\bigr\}.
	\end{align*}
	In general, $\mV_h$ is not a subspace of $\mV$.
	Let $\vQ_h: \vL^2 \rightarrow \mV_h$ denote the $L^2$-orthogonal projection, which is defined by
	\begin{align}
		\bigl(\vv - \vQ_h\vv, \pphi_h\bigr) = 0\qquad\forall \pphi_h\in \mV_h.
	\end{align}
	It is well-known \cite{Girault_Raviart86} that $\vQ_h$ satisfies the following estimates:
	\begin{align}
		\|\vv - \vQ_h\vv\|_{\vL^2} + h\|\nab(\vv - \vQ_h\vv)\|_{\vL^2} &\leq C h^2\|\vA\vv\|_{\vL^2}\qquad\forall \vv\in \mV\cap\vH^2(D),\\
		\|\vv - \vQ_h\vv\|_{\vL^2} &\leq Ch\|\nab\vv\|_{\vL^2}\qquad\forall \vv\in \mV\cap\vH^1(D).
	\end{align} 
	Similarly, let $P_h: L^2_{per}(D) \rightarrow L_h$ be the scalar $L^2$-orthogonal projection, defined by
	\begin{align}
		\bigl(\psi - P_h\psi, q_h\bigr) = 0\qquad\forall q_h \in L_h,
	\end{align}
	then there also holds
	\begin{align}
		\|\psi - P_h\psi\|_{L^2} \leq Ch\|\nab \psi\|_{\vL^2}\qquad\forall \psi \in L^2_{per}(D)\cap H^1(D).
	\end{align}
	Next, we also introduce the well-known trilinear form \cite{Temam}
	\begin{align}
		\tilde{b}(\vu,\vv,\vw) = \bigl(\vu\cdot\nab\vv,\vw\bigr) + \frac12 \bigl([\div \vu]\vv,\vw\bigr)\qquad\forall \vu, \vv, \vw \in \vH^1(D),
	\end{align}
	which is anti-symmetric in the sense that
	\begin{align}
		\tilde{b}(\vu,\vv,\vw) = - \tilde{b}(\vu,\vw,\vv) \qquad\forall\vu,\vv,\vw \in \vH^1(D).
	\end{align}
	Therefore, $$\displaystyle \tilde{b}(\vu, \pphi,\pphi) = 0\qquad\forall\vu,\pphi \in \vH^1(D).$$
	
	Our fully discrete mixed finite method is defined by the following algorithm. 

	\textbf{Algorithm 2.} Let $\vu_h^0$ be a given $\mH_h$-valued random variable. Find $\displaystyle \bigl(\vu_h^{n+1},p_h^{n+1}\bigr) \in L^2(\Ome;\mH_h\times L_h)$ such that there holds $\mP$-a.s.
	\begin{align}
		\label{equu4.5}	&\bigl(\vu_h^{n+1} - \vu_h^n,\pphi_h\bigr) + \nu k \bigl(\nab \vu_h^{n+1},\nab\pphi_h\bigr) + k\bigl(\vu^{n+1}_h\cdot\nab\vu_h^{n+1},\pphi_h\bigr) \\\nonumber
		&\qquad+ \frac{k}{2} \bigl([\div\vu_h^{n+1}]\vu_h^{n+1},\pphi_h\bigr) - k\bigl(p_h^{n+1},\div \pphi_h\bigr) = \bigl(\vg_h^n\Delta W_{n+1},\pphi_h\bigr),\\
		&\qquad\qquad\qquad\qquad\qquad\qquad\qquad\quad\bigl(\div\vu_h^n,\psi_h\bigr) =0,
	\end{align}
	for all $\pphi_h \in \mH_h$ and $\psi_h\in L_h$. $\vg_h^n$ denotes an approximation of $\vg^n$, which satisfies $\|\vg^n - \vg^n_h\|_{\vL^2} \leq Ch^2$.
	
	\smallskip
	Next, we recall the following stability estimates for $\{\vu_h^n\}$ from \cite[Lemma 3.1]{BCP12}. 
	\begin{lemma}\label{stability_FEMs}
		Let $1 \leq q < \infty$ and $\vu_h^0 \in L^{2^q}(\Ome;\mH_h)$ such that $\displaystyle \mE\bigl[\|\vu_h^0\|^{2^q}_{\vL^2}\bigr] \leq C$. Then, there exists a pair $\bigl\{(\vu_h^{n},p^{n}_h)\bigr\}_{n=1}^M \subset L^2(\Ome; \mH_h\times L_h)$ that solves Algorithm 2 and satisfies
		\begin{enumerate}[{\rm (i)}]
			\item $\displaystyle \mE\biggl[\max_{1\leq n \leq M}\|\vu_h^n\|^{2^q}_{\vL^2} + \nu k\sum_{n=1}^M \|\vu_h^n\|^{2^{q}-2}_{\vL^2}\|\nab\vu_h^n\|^2_{\vL^2}\biggr] \leq C_{T,q},$
			\item $\displaystyle \mE\biggl[\biggl(k\sum_{n=1}^M \|\nab\vu_h^n\|^2_{\vL^2}\biggr)^{2^{q-1}}\,\biggr] \leq C_{T,q},$
		\end{enumerate}
		where $C_{T,q} = C_{T,q}(D_T,q,\vu_h^0)$.
	\end{lemma}

	\subsection{High moment and pathwise error estimates for the fully discrete velocity approximation} 
	In this subsection, we state and prove another  main result of this paper.
	
	\begin{theorem}\label{theorem_fully_chapter5} 
		Let $2\leq q < \infty$ and $\vu_h^0 = \vQ_h \vu_0$.  Assume that $\vu_0 \in L^{q}(\Ome;\mV)$. Let $\{(\vu^n, p^n)\}$ and $\{\vu^n_h,p^n_h\}$ be the velocity and pressure approximations generated by Algorithm 1 and 2, respectively. Then there holds
		\begin{align*}
			\bigl(\mE\bigl[\max_{1\leq n \leq M}\|\vu^n - \vu_h^n\|^{q}_{\vL^2}\bigr]\bigr)^{\frac1q} &+ \biggl(\mE\biggl[\Bigl(\nu k \sum_{n=1}^M\|\nab(\vu^n - \vu_h^n)\|^2_{\vL^2}\Bigr)^{\frac{q}{2}}\biggr]\biggr)^{\frac1q}\\\nonumber
			&\leq C_q\biggl( k^{\frac12} + h + k^{-\frac{1}{2}}h\,\biggr),
		\end{align*}
		where $C_q = C(T,q,\vu_0,\vf)>0$ is independent of $k$ and $h$.
	\end{theorem}
	
	
	\bigskip
	
	The proof of Theorem \ref{theorem_fully_chapter5} has several steps and will be accomplished through the following two lemmas.
	\begin{lemma}\label{lemma4.3} Under the assumptions of Theorem \ref{theorem_fully_chapter5}, there holds $\mP$-a.s.
		\begin{align*}
			\max_{1\leq \ell \leq M}&\biggl[\|\vQ_h(\vu^{\ell} - \vu_h^{\ell})\|^2_{\vL^2} + k\sum_{n=1}^{\ell}\|\nab(\vu^{n} - \vu_h^n)\|^2_{\vL^2}\biggr] \leq \bigl(A_h + B_h\bigr)\exp\biggl(k\sum_{n=1}^{M} D^n_h\biggr),
		\end{align*}
		where
		\begin{align}\label{equu4.13}
			A_h &:= \|\vQ_h\vE^0\|^2_{\vL^2} + Ch^2  k\sum_{n=1}^{M}\|\vA\vu^{n}\|^2_{\vL^2} + Ch^3 k\sum_{n=1}^{M}\|\nab\vu^{n}\|^2_{\vL^2}\|\vA\vu^{n}\|^2_{\vL^2}\\\nonumber
			&\qquad+ Ch^2 k\sum_{n=1}^{M}\|\vu_h^{n}\|^2_{\vL^2}\|\vA\vu^{n}\|^2_{\vL^2}   \\\nonumber
			&\qquad+ Ck\max_{1 \leq n \leq M}\|\nab\vu^{n}\|^4_{\vL^2}\Bigl(k\sum_{n = 1}^{M}\|\nab\vu_h^n\|^2_{\vL^2} + k\sum_{n = 1}^M\|\nab\vu^n\|^2_{\vL^2}\Bigr)\\\nonumber
			&\qquad +\frac{16C_L^2k}{\nu}\max_{1 \leq n \leq M}\|\nab\vu^{n}\|^2_{\vL^2}\|\vQ_h\vE^{n}\|^2_{\vL^2} + Ch^2 k \sum_{n=1}^{M}\|\nab p^{n}\|^2_{\vL^2},\\\nonumber
			B_h &:= \sum_{n=1}^{M} \|(\vg^{n-1}- \vg_h^{n-1})\Delta W_{n}\|^2_{\vL^2} \\\nonumber
			&\qquad+  \max_{1\leq \ell \leq M}\Bigl|\sum_{n = 1}^{\ell}\bigl((\vg^{n-1} - \vg_h^{n-1})\Delta W_{n},\vQ_h\vE^{n-1}\bigr)\Bigr|,\\\nonumber
			D_h^n &:= \frac{16C_L^2}{\nu}\|\nab\vu^n\|^2_{\vL^2} + \frac{2\tilde{C}}{\nu}\|\vA\vu^n\|^2_{\vL^2}.
		\end{align}
	\end{lemma}

	\begin{proof} Let $\vE^n := \vu^n - \vu_h^n$ for all $0 \leq n \leq M$. Subtracting \eqref{equu3.1} from \eqref{equu4.5}, we obtain the following error equation
		\begin{align}\label{equu4.14}
			\bigl(\vE^{n+1} &- \vE^{n},\pphi_h\bigr) + \nu k \bigl(\nab\vE^{n+1},\nab\pphi_h\bigr) + k\,\tilde{b}(\vu^{n+1},\vu^{n+1},\pphi_h) \\\nonumber
			&\qquad- k \,\tilde{b}(\vu_h^{n+1},\vu_h^{n+1},\pphi_h)
			- k\bigl(p^{n+1} - p^{n+1}_h,\div\pphi_h\bigr) \\\nonumber
			&\qquad= \bigl((\vg^n - \vg_h^n)\Delta W_{n+1},\pphi_h\bigr)\qquad\forall\pphi_h\in\mH_h.
		\end{align}
		Setting $\pphi_h = \vQ_h \vE^{n+1} \in \mV_h$ in \eqref{equu4.14}, using the orthogonality of the $\vL^2$-projection and the fact that $\bigl(p^{n+1}_h,\div\vQ_h\vE^{n+1}\bigr) = 0$,  we obtain
		\begin{align}\label{equu4.15}
			\bigl(\vQ_h\vE^{n+1} &- \vQ_h\vE^{n},\vQ_h\vE^{n+1}\bigr) + \nu k\|\nab\vE^{n+1}\|^2_{\vL^2}   \\\nonumber
			&= \nu k \bigl(\nab\vE^{n+1}, \nab(\vu^{n+1} - \vQ_h\vu^{n+1})\bigr) - k \,\tilde{b}(\vu^{n+1},\vu^{n+1},\vQ_h\vE^{n+1}) \\\nonumber
			&\qquad+ k \,\tilde{b}(\vu_h^{n+1},\vu_h^{n+1},\vQ_h\vE^{n+1}) + 	 k\bigl(p^{n+1},\div\vQ_h\vE^{n+1}\bigr) \\\nonumber
			&\qquad+ \bigl((\vg^n - \vg_h^n)\Delta W_{n+1},\vQ_h\vE^{n+1}\bigr).
		\end{align}
		We use the identity $2(a,a-b) = \|a\|^2 - \|b\|^2 + \|a-b\|^2$ to rewrite the left-hand side of \eqref{equu4.15} as
		\begin{align}\label{equu4.16}
			\frac12&\bigl[\|\vQ_h\vE^{n+1}\|^2_{\vL^2} - \|\vQ_h\vE^n\|^2_{\vL^2} + \|\vQ_h(\vE^{n+1} - \vE^n)\|^2_{\vL^2}\bigr] + \nu k\|\nab\vE^{n+1}\|^2_{\vL^2}   \\\nonumber
			&= \nu k \bigl(\nab\vE^{n+1}, \nab(\vu^{n+1} - \vQ_h\vu^{n+1})\bigr) -\bigl[ k \,\tilde{b}(\vu^{n+1},\vu^{n+1},\vQ_h\vE^{n+1}) \\\nonumber
			&\qquad - k \,\tilde{b}(\vu_h^{n+1},\vu_h^{n+1},\vQ_h\vE^{n+1})\bigr] + 	 k\bigl(p^{n+1},\div\vQ_h\vE^{n+1}\bigr) \\\nonumber
			&\qquad+ \bigl((\vg^n - \vg_h^n)\Delta W_{n+1},\vQ_h\vE^{n+1}\bigr)\\\nonumber
			&= {\tt I + II + III + IV}.
		\end{align}
		We note that the term {\tt II} needs a special treatment but other terms can be easily controlled as follows:
		\begin{align}
			{\tt I} \leq \frac{\nu k}{16}\|\nab\vE^{n+1}\|^2_{\vL^2} + Ck h^2\|\vA\vu^{n+1}\|^2_{\vL^2}.
		\end{align}
		On noting that $\bigl(P_h p^{n+1}, \div \vQ_h\vE^{n+1}\bigr) = 0$, we have
		\begin{align}
			{\tt III} &= k\bigl(p^{n+1},\div\vQ_h\vE^{n+1}\bigr)\\\nonumber
			&=k\bigl(p^{n+1} - P_h p^{n+1}, \div\vQ_h\vE^{n+1}\bigr)\\\nonumber
			&\leq \frac{\nu k}{16}\|\nab \vE^{n+1}\|^2_{\vL^2} + Ck h^2\|\nab p^{n+1}\|^2_{\vL^2}.
		\end{align}
		\begin{align}
			{\tt IV} &= \bigl((\vg^n - \vg_h^n)\Delta W_{n+1},\vQ_h(\vE^{n+1}-\vE^n)\bigr) + \bigl((\vg^n - \vg_h^n)\Delta W_{n+1},\vQ_h\vE^n\bigr) \\\nonumber
			&\leq 2\|(\vg^n - \vg^n_h)\Delta W_{n+1}\|^2_{\vL^2} + \frac18\|\vQ_h(\vE^{n+1} -\vE^n)\|^2_{\vL^2} \\\nonumber
			&\qquad+ \bigl((\vg^n - \vg_h^n)\Delta W_{n+1},\vQ_h\vE^n\bigr)
		\end{align}
		To bound the term {\tt II}, we first rewrite it as
		\begin{align}\label{eqII}
			{\tt II} &= -k\bigl[ \tilde{b}(\vu^{n+1},\vu^{n+1},\vQ_h\vE^{n+1}) - \tilde{b}(\vu^{n+1}_h, \vu_h^{n+1},\vQ_h\vE^{n+1})\bigr]\\\nonumber
			&= -k\tilde{b}(\vu^{n+1},\vE^{n+1},\vQ_h\vE^{n+1}) - k\tilde{b}(\vE^{n+1},\vu_h^{n+1},\vQ_h\vE^{n+1})\\\nonumber
			&= -k\tilde{b}(\vu^{n+1},\vE^{n+1},\vE^{n+1}) + k\tilde{b}(\vu^{n+1}, \vE^{n+1}, \vu^{n+1} - \vQ_h \vu^{n+1}) \\\nonumber
			&\qquad+ k\tilde{b}(\vE^{n+1},\vE^{n+1},\vQ_h\vE^{n+1}) - k\tilde{b}(\vE^{n+1},\vu^{n+1},\vQ_h\vE^{n+1})\\\nonumber
			&= -k\tilde{b}(\vu^{n+1},\vE^{n+1},\vE^{n+1}) + k\tilde{b}(\vu^{n+1}, \vE^{n+1}, \vu^{n+1} - \vQ_h \vu^{n+1}) \\\nonumber
			&\qquad+ k\tilde{b}(\vE^{n+1},\vE^{n+1},\vE^{n+1}) - k\tilde{b}(\vE^{n+1},\vE^{n+1},\vu^{n+1} - \vQ_h\vu^{n+1}) \\\nonumber
			&\qquad\qquad- k\tilde{b}(\vE^{n+1},\vu^{n+1},\vQ_h\vE^{n+1})\\\nonumber
			&= k\tilde{b}(\vu^{n+1}, \vE^{n+1}, \vu^{n+1} - \vQ_h \vu^{n+1}) - k\tilde{b}(\vE^{n+1},\vE^{n+1},\vu^{n+1} - \vQ_h\vu^{n+1}) \\\nonumber
			&\qquad\qquad- k\tilde{b}(\vE^{n+1},\vu^{n+1},\vQ_h\vE^{n+1})\\\nonumber
			& := {\tt II_1 + II_2 + II_3}.
		\end{align}
		Here, in order to obtain the last line of \eqref{eqII}, we have used the facts that
		$\tilde{b}(\vu^{n+1},\vE^{n+1},\vE^{n+1}) = 0$ and $\tilde{b}(\vE^{n+1},\vE^{n+1},\vE^{n+1}) = 0$.
		Then we have
		\begin{align*}
			{\tt II_1} &= k\,\tilde{b}(\vu^{n+1},\vE^{n+1}, \vu^{n+1} - \vQ_h\vu^{n+1})\\\nonumber
			&= k\bigl(\vu^{n+1}\cdot\nab\vE^{n+1},\vu^{n+1} - \vQ_h\vu^{n+1}\bigr)+ \frac12 k\bigl([\div\vu^{n+1}]\vE^{n+1},\vu^{n+1} - \vQ_h\vu^{n+1}\bigr)\\\nonumber
			&\leq k\|\vu^{n+1}\|_{\vL^4}\|\nab\vE^{n+1}\|_{\vL^2}\|\vu^{n+1} - \vQ_h\vu^{n+1}\|_{\vL^4}\\\nonumber
			&\leq \frac{\nu k}{16} \|\nab\vE^{n+1}\|^2_{\vL^2} + Ck\|\vu^{n+1}\|^2_{\vL^4}\|\vu^{n+1} - \vQ_h\vu^{n+1}\|^2_{\vL^4}\\\nonumber
			&\leq \frac{\nu k}{16} \|\nab\vE^{n+1}\|^2_{\vL^2} + Ckh^3\|\nab\vu^{n+1}\|^2_{\vL^2}\|\vA\vu^{n+1}\|^2_{\vL^2}. \\
			{\tt II_2} &= -k\,\tilde{b}(\vE^{n+1},\vE^{n+1},\vu^{n+1} - \vQ_h\vu^{n+1})\\\nonumber
			&= -k\bigl(\vE^{n+1}\cdot\nab\vE^{n+1},\vu^{n+1} - \vQ_h\vu^{n+1}\bigr) -\frac{k}{2}\bigl([\div \vE^{n+1}]\vE^{n+1},\vu^{n+1} - \vQ_h\vu^{n+1}\bigr)\\\nonumber
			&\leq 2k\|\vE^{n+1}\|_{\vL^4}\|\nab\vE^{n+1}\|_{\vL^2}\|\vu^{n+1} - \vQ_h\vu^{n+1}\|_{\vL^4}\\\nonumber
			&\leq \frac{\nu k}{16}\|\nab\vE^{n+1}\|^2_{\vL^2} + Ck\|\vE^{n+1}\|^2_{\vL^4}\|\vu^{n+1} - \vQ_h\vu^{n+1}\|^2_{\vL^4}\\\nonumber
			&\leq \frac{\nu k}{16}\|\nab\vE^{n+1}\|^2_{\vL^2} + Ckh^3\|\vE^{n+1}\|^2_{\vL^4}\|\vA\vu^{n+1}\|^2_{\vL^2}\\\nonumber
			&\leq \frac{\nu k}{16}\|\nab\vE^{n+1}\|^2_{\vL^2} + Ckh^3\|\vu^{n+1}\|^2_{\vL^4}\|\vA\vu^{n+1}\|^2_{\vL^2} + Ckh^3\|\vu_h^{n+1}\|^2_{\vL^4}\|\vA\vu^{n+1}\|^2_{\vL^2}\\\nonumber
			&\leq \frac{\nu k}{16}\|\nab\vE^{n+1}\|^2_{\vL^2} + Ckh^3\|\nab\vu^{n+1}\|^2_{\vL^2}\|\vA\vu^{n+1}\|^2_{\vL^2} + Ckh^2\|\vu_h^{n+1}\|^2_{\vL^2}\|\vA\vu^{n+1}\|^2_{\vL^2}.
		\end{align*}
		Here, we have used the inverse inequality
		$\|\vu_h^{n+1}\|_{\vL^4} \leq Ch^{-\frac12}\|\vu_h^{n+1}\|_{\vL^2}$ to get the last inequality.
		\begin{align}
			{\tt II_3} &= -k\,\tilde{b}(\vE^{n+1},\vu^{n+1},\vQ_h\vE^{n+1})\\\nonumber
			&=-k\bigl(\vE^{n+1}\cdot\nab\vu^{n+1},\vQ_h\vE^{n+1}\bigr) - \frac{k}{2}\bigl([\div\vE^{n+1}]\vu^{n+1},\vQ_h\vE^{n+1}\bigr)\\\nonumber
			&= -k\bigl(\vE^{n+1}\cdot\nab\vu^{n+1},\vQ_h\vE^{n+1}\bigr) - \frac{k}{2}\bigl([\div\vE^{n+1}]\vu^{n+1},\vQ_h(\vE^{n+1} - \vE^{n})\bigr) \\\nonumber
			&\qquad+ \frac{k}{2}\bigl([\div\vE^{n+1}]\vu^{n+1},\vQ_h\vE^n\bigr) \\\nonumber
			&={\tt II_{3a} + II_{3b} + II_{3c}}.
		\end{align}
		We can bound each term on the right-hand side of ${\tt II_3}$ as follows.
		\begin{align*}
			{\tt II_{3a}} &= -k\bigl(\vQ_h\vE^{n+1}\cdot\nab\vu^{n+1},\vQ_h\vE^{n+1}\bigr)- k\bigl((\vu^{n+1}-\vQ_h\vu^{n+1})\cdot\nab\vu^{n+1},\vQ_h\vE^{n+1}\bigr)\\\nonumber
			&\leq k\|\vQ_h\vE^{n+1}\|^2_{\vL^4}\|\nab\vu^{n+1}\|_{\vL^2} + k\|\vu^{n+1}-\vQ_h\vu^{n+1}\|_{\vL^4}\|\nab\vu^{n+1}\|_{\vL^2}\|\vQ_h\vE^{n+1}\|_{\vL^4}\\\nonumber
			&\leq C_Lk\|\vQ_h\vE^{n+1}\|_{\vL^2}\|\nab\vE^{n+1}\|_{\vL^2}\|\nab\vu^{n+1}\|_{\vL^2}\\\nonumber
			&\qquad+Ck\|\vu^{n+1}-\vQ_h\vu^{n+1}\|_{\vL^4}\|\nab\vu^{n+1}\|_{\vL^2}\|\nab\vE^{n+1}\|_{\vL^2}\\\nonumber
			&\leq \frac{8C_L^2}{\nu}k\|\nab\vu^{n+1}\|^2_{\vL^2}\|\vQ_h\vE^{n+1}\|^2_{\vL^2} \\\nonumber
			&\qquad+ Ck\|\vu^{n+1} - \vQ_h\vu^{n+1}\|^2_{\vL^4}\|\nab\vu^{n+1}\|^2_{\vL^2} + \frac{\nu k}{16}\|\nab\vE^{n+1}\|^2_{\vL^2}\\\nonumber
			&\leq \frac{8C_L^2}{\nu} k\|\nab\vu^{n+1}\|^2_{\vL^2}\|\vQ_h\vE^{n+1}\|^2_{\vL^2} + Ckh^3\|\nab\vu^{n+1}\|^2_{\vL^2}\|\vA\vu^{n+1}\|^2_{\vL^2} + \frac{\nu k}{16}\|\nab\vE^{n+1}\|^2_{\vL^2}.\\
			{\tt II_{3b}} &= -\frac{k}{2}\bigl([\div\vE^{n+1}]\vu^{n+1},\vQ_h(\vE^{n+1} - \vE^n)\bigr)\\\nonumber
			&\leq \frac{k}{2}\|\nab\vE^{n+1}\|_{\vL^2}\|\vu^{n+1}\|_{\vL^4}\|\vQ_h(\vE^{n+1} - \vE^n)\|_{\vL^4}\\\nonumber
			&\leq C k \|\nab\vE^{n+1}\|_{\vL^2}\|\nab\vu^{n+1}\|_{\vL^2}\|\vQ_h(\vE^{n+1} - \vE^n)\|^{1/2}_{\vL^2}\|\nab (\vE^{n+1} - \vE^n)\|^{1/2}_{\vL^2}\\\nonumber
			&\leq \frac{\nu k}{16}\|\nab\vE^{n+1}\|^2_{\vL^2} + Ck\|\nab\vu^{n+1}\|^2_{\vL^2}\|\vQ_h(\vE^{n+1} - \vE^n)\|_{\vL^2}\|\nab (\vE^{n+1} - \vE^n)\|_{\vL^2}\\\nonumber
			&\leq \frac{\nu k}{16}\|\nab\vE^{n+1}\|^2_{\vL^2} + \frac{1}{8}\|\vQ_h(\vE^{n+1} - \vE^n)\|^2_{\vL^2} \\\nonumber
			&\qquad+ Ck^2\|\nab\vu^{n+1}\|^4_{\vL^2}\|\nab(\vE^{n+1} - \vE^n)\|^2_{\vL^2}.\\\nonumber
			{\tt II_{3c}} &= \frac{k}{2}\bigl([\div\vE^{n+1}]\vu^{n+1},\vQ_h\vE^n\bigr)\\\nonumber
			&\leq \frac{k}{2}\|\nab\vE^{n+1}\|_{\vL^2}\|\vu^{n+1}\cdot \vQ_h\vE^n\|_{\vL^2}\\\nonumber
			&\leq \frac{k}{2}\|\nab\vE^{n+1}\|_{\vL^2}\|\vu^{n+1}\|_{\vL^{\infty}}\|\vQ_h\vE^n\|_{\vL^2}\\\nonumber
			&\leq \frac{\nu k}{16}\|\nab\vE^{n+1}\|^2_{\vL^2} + \frac{k}{\nu}\|\vu^{n+1}\|^2_{\vL^{\infty}}\|\vQ_h\vE^n\|^2_{\vL^2}\\\nonumber
			&\leq \frac{\nu k}{16}\|\nab\vE^{n+1}\|^2_{\vL^2} + \frac{\tilde{C}}{\nu} k\|\vA\vu^{n+1}\|^2_{\vL^2}\|\vQ_h\vE^n\|^2_{\vL^2},
		\end{align*}
		where the last inequality was obtained by using the Sobolev embedding: $\|\vu^{n+1}\|^2_{\vL^{\infty}} \leq \tilde{C}\|\vA\vu^{n+1}\|^2_{\vL^2}$ for some constant $\tilde{C}>0$.
		
		Now, substituting the above estimates for the terms {\tt I}--{\tt IV} into the right-hand side of  the right-hand side of \eqref{equu4.16}, we obtain
		\begin{align}\label{equu4.31}
			\frac12\bigl[&\|\vQ_h\vE^{n+1}\|^2_{\vL^2} - \|\vQ_h\vE^n\|^2_{\vL^2}\bigr] +\frac14\|\vQ_h(\vE^{n+1} - \vQ_h\vE^n)\|^2_{\vL^2} \\\nonumber
			&\qquad\qquad\qquad\qquad\qquad\qquad\qquad+ \frac{9\nu k}{16}\|\nab\vE^{n+1}\|^2_{\vL^2}\\\nonumber
			&\leq Ckh^2\|\vA\vu^{n+1}\|^2_{\vL^2} + Ckh^3\|\nab\vu^{n+1}\|^2_{\vL^2}\|\vA\vu^{n+1}\|^2_{\vL^2}\\\nonumber
			&\qquad+ Ckh^2\|\vu_h^{n+1}\|^2_{\vL^2}\|\vA\vu^{n+1}\|^2_{\vL^2} + \frac{8C_L^2}{\nu} k\|\nab\vu^{n+1}\|^2_{\vL^2}\|\vQ_h\vE^{n+1}\|^2_{\vL^2}\\\nonumber
			&\qquad+ Ck^2\|\nab\vu^{n+1}\|^4_{\vL^2}\|\nab(\vE^{n+1} - \vE^n)\|^2_{\vL^2} \\\nonumber
			&\qquad+ \frac{\tilde{C}k}{\nu}\|\vA\vu^{n+1}\|^2_{\vL^2}\|\vQ_h\vE^n\|^2_{\vL^2} + Ck h^2\|\nab p^{n+1}\|^2_{\vL^2}\\\nonumber
			&\qquad+ \|(\vg^n - \vg^n_h)\Delta W_{n+1}\|^2_{\vL^2} + \bigl((\vg^n - \vg_h^n)\Delta W_{n+1},\vQ_h\vE^n\bigr).
		\end{align} 
		
		Lowering the index $n$ in \eqref{equu4.31} by $1$ and taking the summation operator $\displaystyle \sum_{n=1}^{\ell}$ on both sides yield
		\begin{align}\label{equu4.32}
			\frac12&\|\vQ_h\vE^{\ell}\|^2_{\vL^2} + \frac14\sum_{n=1}^{\ell}\|\vQ_h(\vE^{n} - \vE^{n-1})\|^2_{\vL^2} + \frac{9\nu k}{16} \sum_{n=1}^{\ell}\|\nab\vE^{n}\|^2_{\vL^2} \\\nonumber
			&\leq \frac12\|\vQ_h\vE^0\|^2_{\vL^2} + Ch^2  k\sum_{n=1}^{\ell}\|\vA\vu^{n}\|^2_{\vL^2} + Ch^3 k\sum_{n=1}^{\ell}\|\nab\vu^{n}\|^2_{\vL^2}\|\vA\vu^{n}\|^2_{\vL^2}\\\nonumber
			&\qquad+ Ch^2 k\sum_{n=1}^{\ell}\|\vu_h^{n}\|^2_{\vL^2}\|\vA\vu^{n}\|^2_{\vL^2} + \frac{8C_L^2}{\nu} k\sum_{n=1}^{\ell}\|\nab\vu^{n}\|^2_{\vL^2}\|\vQ_h\vE^{n}\|^2_{\vL^2} \\\nonumber
			&\qquad+ Ck^2 \sum_{n=1}^{\ell}\|\nab\vu^{n}\|^4_{\vL^2}\bigl(\|\nab(\vu_h^{n} - \vu_h^{n-1})\|^2_{\vL^2} + \|\nab(\vu^n - \vu^{n-1})\|^2_{\vL^2}\bigr)\\\nonumber
			&\qquad+ \frac{\tilde{C}}{\nu} k  \sum_{n=1}^{\ell} \|\vA\vu^{n}\|^2_{\vL^2}\|\vQ_h\vE^{n-1}\|^2_{\vL^2} + Ch^2 k \sum_{n=1}^{\ell}\|\nab p^{n}\|^2_{\vL^2}\\\nonumber
			&\qquad+ \sum_{n=1}^{\ell} \|(\vg^{n-1}- \vg_h^{n-1})\Delta W_{n}\|_{\vL^2}^2 + \sum_{n = 1}^{\ell}\bigl((\vg^{n-1} - \vg_h^{n-1})\Delta W_{n},\vQ_h\vE^{n-1}\bigr).
		\end{align}
		
		Bounding the sixth term on the right-hand side of \eqref{equu4.32} as follows
		\begin{align*}
			&Ck^2 \sum_{n=1}^{\ell}\|\nab\vu^{n}\|^4_{\vL^2}\bigl(\|\nab(\vu_h^{n} - \vu_h^{n-1})\|^2_{\vL^2} + \|\nab(\vu^n - \vu^{n-1})\|^2_{\vL^2}\bigr)\\\nonumber
			&\leq Ck\max_{1 \leq n \leq M}\|\nab\vu^{n}\|^4_{\vL^2}\Bigl(k\sum_{n = 1}^{M}\|\nab\vu_h^n\|^2_{\vL^2} + k\sum_{n = 1}^M\|\nab\vu^n\|^2_{\vL^2}\Bigr),
		\end{align*}
		then, \eqref{equu4.32} becomes
		\begin{align}\label{equu4.322}
			\frac12&\|\vQ_h\vE^{\ell}\|^2_{\vL^2} + \frac14\sum_{n=1}^{\ell}\|\vQ_h(\vE^{n} - \vE^{n-1})\|^2_{\vL^2} + \frac{9\nu k}{16} \sum_{n=1}^{\ell}\|\nab\vE^{n}\|^2_{\vL^2} \\\nonumber
			&\leq \frac12\|\vQ_h\vE^0\|^2_{\vL^2} + Ch^2  k\sum_{n=1}^{\ell}\|\vA\vu^{n}\|^2_{\vL^2} + Ch^3 k\sum_{n=1}^{\ell}\|\nab\vu^{n}\|^2_{\vL^2}\|\vA\vu^{n}\|^2_{\vL^2}\\\nonumber
			&\qquad+ Ch^2 k\sum_{n=1}^{\ell}\|\vu_h^{n}\|^2_{\vL^2}\|\vA\vu^{n}\|^2_{\vL^2} + \frac{8C_L^2}{\nu} k\sum_{n=1}^{\ell}\|\nab\vu^{n}\|^2_{\vL^2}\|\vQ_h\vE^{n}\|^2_{\vL^2} \\\nonumber
			&\qquad+ Ck\max_{1 \leq n \leq M}\|\nab\vu^{n}\|^4_{\vL^2}\Bigl(k\sum_{n = 1}^{M}\|\nab\vu_h^n\|^2_{\vL^2} + k\sum_{n = 1}^M\|\nab\vu^n\|^2_{\vL^2}\Bigr)\\\nonumber
			&\qquad+ \frac{\tilde{C}}{\nu} k  \sum_{n=1}^{\ell} \|\vA\vu^{n}\|^2_{\vL^2}\|\vQ_h\vE^{n-1}\|^2_{\vL^2} + Ch^2 k \sum_{n=1}^{\ell}\|\nab p^{n}\|^2_{\vL^2}\\\nonumber
			&\qquad+ \sum_{n=1}^{\ell} \|(\vg^{n-1}- \vg_h^{n-1})\Delta W_{n}\|_{\vL^2}^2 + \sum_{n = 1}^{\ell}\bigl((\vg^{n-1} - \vg_h^{n-1})\Delta W_{n},\vQ_h\vE^{n-1}\bigr)\\\nonumber
			&\leq \frac12\|\vQ_h\vE^0\|^2_{\vL^2} + Ch^2  k\sum_{n=1}^{\ell}\|\vA\vu^{n}\|^2_{\vL^2} + Ch^3 k\sum_{n=1}^{\ell}\|\nab\vu^{n}\|^2_{\vL^2}\|\vA\vu^{n}\|^2_{\vL^2}\\\nonumber
			&\qquad+ Ch^2 k\sum_{n=1}^{\ell}\|\vu_h^{n}\|^2_{\vL^2}\|\vA\vu^{n}\|^2_{\vL^2} \\\nonumber
			&\qquad+  k\sum_{n=1}^{\ell-1}\Bigl(\frac{8C_L^2}{\nu}\|\nab\vu^{n}\|^2_{\vL^2} + \frac{\tilde{C}}{\nu}\|\vA\vu^{n+1}\|^2_{\vL^2}\Bigr)\|\vQ_h\vE^{n}\|^2_{\vL^2} \\\nonumber
			&\qquad+ Ck\max_{1 \leq n \leq M}\|\nab\vu^{n}\|^4_{\vL^2}\Bigl(k\sum_{n = 1}^{M}\|\nab\vu_h^n\|^2_{\vL^2} + k\sum_{n = 1}^M\|\nab\vu^n\|^2_{\vL^2}\Bigr)\\\nonumber
			&\qquad+ \frac{8C_L^2k}{\nu} \max_{1 \leq n \leq M}\|\nab\vu^n\|^2_{\vL^2}\|\vQ_h\vE^n\|^2_{\vL^2} + Ch^2 k \sum_{n=1}^{\ell}\|\nab p^{n}\|^2_{\vL^2}\\\nonumber
			&\qquad+ \sum_{n=1}^{\ell} \|(\vg^{n-1}- \vg_h^{n-1})\Delta W_{n}\|_{\vL^2}^2 + \sum_{n = 1}^{\ell}\bigl((\vg^{n-1} - \vg_h^{n-1})\Delta W_{n},\vQ_h\vE^{n-1}\bigr)\\\nonumber
			&\leq \frac12\|\vQ_h\vE^0\|^2_{\vL^2} + Ch^2  k\sum_{n=1}^{M}\|\vA\vu^{n}\|^2_{\vL^2} + Ch^3 k\sum_{n=1}^{M}\|\nab\vu^{n}\|^2_{\vL^2}\|\vA\vu^{n}\|^2_{\vL^2}\\\nonumber
			&\qquad+ Ch^2 k\sum_{n=1}^{M}\|\vu_h^{n}\|^2_{\vL^2}\|\vA\vu^{n}\|^2_{\vL^2} \\\nonumber
			&\qquad+  k\sum_{n=1}^{\ell-1}\Bigl(\frac{8C_L^2}{\nu}\|\nab\vu^{n}\|^2_{\vL^2} + \frac{\tilde{C}}{\nu}\|\vA\vu^{n+1}\|^2_{\vL^2}\Bigr)\|\vQ_h\vE^{n}\|^2_{\vL^2} \\\nonumber
			&\qquad+ Ck\max_{1 \leq n \leq M}\|\nab\vu^{n}\|^4_{\vL^2}\Bigl(k\sum_{n = 1}^{M}\|\nab\vu_h^n\|^2_{\vL^2} + k\sum_{n = 1}^M\|\nab\vu^n\|^2_{\vL^2}\Bigr)\\\nonumber
			&\qquad+ \frac{8C_L^2k}{\nu} \max_{1 \leq n \leq M}\|\nab\vu^n\|^2_{\vL^2}\|\vQ_h\vE^n\|^2_{\vL^2} + Ch^2 k \sum_{n=1}^{M}\|\nab p^{n}\|^2_{\vL^2}\\\nonumber
			&\qquad+ \sum_{n=1}^{M} \|(\vg^{n-1}- \vg_h^{n-1})\Delta W_{n}\|_{\vL^2}^2 \\\nonumber
			&\qquad+ \max_{1\leq \ell \leq M}\Bigl|\sum_{n = 1}^{\ell}\bigl((\vg^{n-1} - \vg_h^{n-1})\Delta W_{n},\vQ_h\vE^{n-1}\bigr)\Bigr|,
		\end{align}
		which is equivalent to (see the notations defined in \eqref{equu4.16})
		\begin{align}\label{equu4.34}
			&\|\vQ_h\vE^{\ell}\|^2_{\vL^2} + \sum_{n=1}^{\ell}\|\vQ_h(\vE^{n} - \vE^{n-1})\|^2_{\vL^2} + \nu k \sum_{n=1}^{\ell}\|\nab\vE^{n}\|^2_{\vL^2}\\\nonumber
			&\leq A_h +B_h + k\sum_{n=1}^{\ell-1}D^n_h\max_{1 \leq \ell \leq n}\|\vQ_h\vE^{\ell}\|^2_{\vL^2}\\\nonumber
			& \leq \bigl(A_h + B_h\bigr) \exp\Bigl( k\sum_{n=1}^{m-1}D^n_h\Bigr).
		\end{align} 
		
		The proof is completed by replacing $m$ by $M$ in \eqref{equu4.34}.
	\end{proof}
	
	\smallskip
	\begin{lemma}\label{lemma4.4} Under the assumptions of Theorem \ref{theorem_fully_chapter5} and Lemma \ref{lemma4.3}, there holds
		\begin{align}
			\biggl(\mE\biggl[\bigl(A_h +& B_h\bigr)^{\frac{q}{2}}\exp\Bigl(\frac{qk}{2}\sum_{n=1}^M D^n_h\Bigr)\biggr]\biggr)^{\frac1q}  \leq C_q \biggl(k^{\frac12} + h + k^{-\frac{1}{2}}h\,\biggr).
		\end{align}	
	\end{lemma}
	\begin{proof}
		First, by using Cauchy-Schwarz inequality, we have:
		\begin{align}\label{equu4.35}
			\mE\biggl[\bigl(A_h +& B_h\bigr)^{\frac{q}{2}}\exp\Bigl( \frac{qk}{2}\sum_{n=1}^M D_h^n\Bigr)\biggr] \\\nonumber
			&\leq C_q\bigl(\mE\bigl[A_h^q\bigr]\bigr)^{\frac{1}{2}}\biggl(\mE\Bigl[\exp\Bigl(qk\sum_{n = 1}^MD^n_h\Bigr)\Bigr]\biggr)^{\frac12} \\\nonumber
			&\qquad+ C_q\bigl(\mE\bigl[B_h^q\bigr]\bigr)^{\frac{1}{2}}\biggl(\mE\Bigl[\exp\Bigl(qk\sum_{n = 1}^MD^n_h\Bigr)\Bigr]\biggr)^{\frac12}.
		\end{align}
		By using Lemma \ref{lemma_exp_discrete}, we can easily control the factor $\biggl(\mE\Bigl[\exp\Bigl(q k\sum_{n=1}^M D^n_h\Bigr)\Bigr]\biggr)^{\frac12}$. In fact, 
		\begin{align*}
			\mE\Bigl[\exp\Bigl(qk\sum_{n=1}^M D^n_h\Bigr)\Bigr] &= \mE\Bigl[\exp\Bigl( k\sum_{n=1}^M \frac{32qC_L^2}{\nu}\|\nab\vu^{n}\|^2_{\vL^2} + k\sum_{n = 1}^M\frac{4q\tilde{C}}{\nu}\|\vA\vu^{n+1}\|^2_{\vL^2}\Bigr)\Bigr] \\\nonumber
			&\leq \mE\Bigl[\exp\Bigl( \frac{32qC_L^2T}{\nu}\max_{1 \leq n \leq M}\|\nab\vu^{n}\|^2_{\vL^2} + k\sum_{n = 1}^M\frac{4q\tilde{C}}{\nu}\|\vA\vu^{n+1}\|^2_{\vL^2}\Bigr)\Bigr] \\\nonumber
			&\leq \mE\Bigl[\exp\Bigl( \sigma\max_{1 \leq n \leq M}\|\nab\vu^{n}\|^2_{\vL^2} + \sigma k\sum_{n = 1}^M\|\vA\vu^{n+1}\|^2_{\vL^2}\Bigr)\Bigr] \leq C_2,
		\end{align*}
		where $\sigma_q = \max\bigl\{\frac{32qC_L^2 T}{\nu}, \frac{4q\tilde{C}}{\nu}\bigr\}$.
		
		Therefore, it remains to bound $\mE[A_h^q]$ and $\mE[B_h^q]$. By the assumptions on $\vg$, we have
		\begin{align}\label{equu4.36}
			\mE[B_h^q] &= \mE\Bigl[\Bigl(\sum_{n=1}^{M} \|(\vg^{n-1} - \vg_h^{n-1})\Delta W_{n}\|^2_{\vL^2} \\\nonumber
			&\qquad+  \max_{1\leq \ell \leq M}\Bigl|\sum_{n = 1}^{\ell}\bigl((\vg^n - \vg_h^n)\Delta W_{n+1},\vQ_h\vE^n\bigr)\Bigr|\Bigr)^q\Bigr]\\\nonumber
			&\leq C_q\mE\Bigl[\Bigl(\sum_{n = 1}^M\|\vg^{n-1}-\vg_h^{n-1}\|^2_{\vL^2}|\Delta W_n|^2\Bigr)^q\Bigr]\\\nonumber
			&\qquad+ C_q\mE\Big[\max_{1\leq \ell \leq M}\Bigl|\sum_{n = 1}^{\ell}\bigl((\vg^n - \vg_h^n)\Delta W_{n+1},\vQ_h\vE^n\bigr)\Bigr|^q\Big]\\\nonumber
			&\leq C_qh^{2q}\mE\Bigl[\Bigl(\sum_{n = 1}^M|\Delta W_n|^2\Bigr)^q\Bigr] \\\nonumber
			&\qquad+ C_q\mE\Big[\max_{1\leq \ell \leq M}\Bigl|\sum_{n = 1}^{\ell}\bigl((\vg^{n-1} - \vg_h^{n-1})\Delta W_{n},\vQ_h\vE^{n-1}\bigr)\Bigr|^q\Big]\\\nonumber
			&= B_{h,1} + B_{h,2}.
		\end{align}
		It follows from the discrete H\"older inequality and \eqref{mean_wiener} that
		\begin{align}
			B_{h,1} &= C_qh^{2q}\mE\Bigl[\Bigl(\sum_{n = 1}^M|\Delta W_n|^2\Bigr)^q\Bigr]\\\nonumber
			&\leq C_qh^{2q}\, M^{q-1}\sum_{n = 1}^M\mE\bigl[|\Delta W_n|^{2q}\bigr]\leq C_qh^{2q} M^{q} k^q \leq C_qh^{2q}.
		\end{align}
		To estimate $B_{h,2}$, we use the Burkholder--Davis--Gundy inequality and the assumptions on $\vg$ to obtain
		\begin{align}
			B_{h,2} &= C_q\mE\Big[\max_{1\leq \ell \leq M}\Bigl|\sum_{n = 1}^{\ell}\bigl((\vg^n - \vg_h^n)\Delta W_{n},\vQ_h\vE^{n-1}\bigr)\Bigr|^q\Big]\\\nonumber
			&\leq C_q\mE\Bigl[\Bigl(k\sum_{n = 1}^M\|\vg^{n-1} - \vg^{n-1}_h\|^2_{\vL^2}\|\vQ_h\vE^{n-1}\|^2_{\vL^2}\Bigr)^{\frac{q}{2}} \Bigr]\\\nonumber
			&\leq C_q h^{2q} \mE\Bigl[\Bigl(k\sum_{n = 1}^M\|\vE^{n-1}\|^2_{\vL^2}\Bigr)^{\frac{q}{2}}\Bigr] \leq C_qh^{2q}.
		\end{align}
		In addition, using the stability estimates from Lemmas \ref{stability_means}, \ref{stability_FEMs} and \ref{stability_pressures} and the assumption $\mE[\|\vQ_h\vE^0\|^q_{\vL^2}] \leq Ch^q$, we obtain
		\begin{align}\label{equu4.37}
			\mE[A_h^q] &= \mE\bigg[\biggl(\|\vQ_h\vE^0\|^2_{\vL^2} + Ch^2  k\sum_{n=1}^{M}\|\vA\vu^{n}\|^2_{\vL^2} + Ch^3 k\sum_{n=1}^{M}\|\nab\vu^{n}\|^2_{\vL^2}\|\vA\vu^{n}\|^2_{\vL^2}\\\nonumber
			&\qquad+ Ch^2 k\sum_{n=1}^{M}\|\vu_h^{n}\|^2_{\vL^2}\|\vA\vu^{n}\|^2_{\vL^2}   \\\nonumber
			&\qquad+ Ck\max_{1 \leq n \leq M}\|\nab\vu^{n}\|^4_{\vL^2}\Bigl(k\sum_{n = 1}^{M}\|\nab\vu_h^n\|^2_{\vL^2} + k\sum_{n = 1}^M\|\nab\vu^n\|^2_{\vL^2}\Bigr)\\\nonumber
			&\qquad +\frac{16C_L^2k}{\nu}\max_{1 \leq n \leq M}\|\nab\vu^{n}\|^2_{\vL^2}\|\vQ_h\vE^{n}\|^2_{\vL^2} + Ch^2 k \sum_{n=1}^{M}\|\nab p^{n}\|^2_{\vL^2}\biggr)^q\bigg]\\\nonumber
			&\leq C_qh^{2q} + C_qh^{2q}\mE\Bigl[\Bigl(k\sum_{n=1}^{M}\|\vA\vu^{n}\|^2_{\vL^2}\Bigr)^q\Bigr]\\\nonumber
			&\qquad+ Ch^{3q}\mE\Bigl[\Bigl(k\sum_{n=1}^{M}\|\nab\vu^{n}\|^2_{\vL^2}\|\vA\vu^{n}\|^{2}_{\vL^2}\Bigr)^q\Bigr] \\\nonumber
			&\qquad+C_qh^{2q}\mE\Bigl[\max_{1 \leq n \leq M}\|\vu_h^{n}\|^{2q}_{\vL^2}\Bigl(k\sum_{n=1}^{M}\|\vA\vu^{n}\|^2_{\vL^2}\Bigr)^{q}\Bigr]\\\nonumber
			&\qquad+C_qk^q\mE\biggl[\max_{1 \leq n \leq M}\|\nab\vu^n\|^{4q}_{\vL^2}\Bigl(k\sum_{n = 1}^M\|\nab\vu^n_h\|^2_{\vL^2}\Bigr)^q\biggr] \\\nonumber
			&\qquad+  C_qk^q\mE\biggl[\max_{1 \leq n \leq M}\|\nab\vu^n\|^{4q}_{\vL^2}\Bigl(k\sum_{n = 1}^M\|\nab\vu^n\|^2_{\vL^2}\Bigr)^{q}\biggr]\\\nonumber
			&\qquad + C_qk^q\mE\Bigl[\max_{1 \leq n \leq M}\|\nab\vu^n\|^{2q}_{\vL^2}\|\vE^n\|^{2q}_{\vL^2}\Bigr] + C_qh^{2q}\mE\Bigl[\Bigl(k \sum_{n=1}^{M}\|\nab p^{n}\|^2_{\vL^2}\Bigr)^{q}\Bigr]\\\nonumber
			&\leq C_qh^{2q} + C_qh^{2q} \Bigl(\mE\Bigl[\max_{1 \leq n \leq M}\|\vu_h^n\|^{4q}_{\vL^2}\Bigr]\Bigr)^{\frac12}\Bigl(\mE\Bigl[\Bigl(k\sum_{n = 1}^M\|\vA\vu^n\|^2_{\vL^2}\Bigr)^{2q}\Bigr]\Bigr)^{\frac12}\\\nonumber
			&\qquad+ C_qk^q \Bigl(\mE\Bigl[\max_{1 \leq n \leq M}\|\nab\vu^n\|^{8q}_{\vL^2}\Bigr]\Bigr)^{\frac12}\Bigl(\mE\Bigl[\Bigl(k\sum_{n = 1}^M\|\nab\vu_h^n\|^2_{\vL^2}\Bigr)^{2q}\Bigr]\Bigr)^{\frac12}\\\nonumber
			&\qquad+ C_qk^q \Bigl(\mE\Bigl[\max_{1 \leq n \leq M}\|\nab\vu^n\|^{8q}_{\vL^2}\Bigr]\Bigr)^{\frac12}\Bigl(\mE\Bigl[\Bigl(k\sum_{n = 1}^M\|\nab\vu^n\|^2_{\vL^2}\Bigr)^{2q}\Bigr]\Bigr)^{\frac12}\\\nonumber
			&\qquad+ C_qk^q \mE\Bigl[\max_{1 \leq n \leq M}\|\nab\vu^n\|^{2q}_{\vL^2}\bigl(\|\vu^n\|^{2q}_{\vL^2} + \|\vu_h^n\|^{2q}_{\vL^2}\bigr)\Bigr]\\\nonumber
			&\qquad+ C_qh^{2q}\mE\Bigl[\Bigl(k \sum_{n=1}^{M}\|\nab p^{n}\|^2_{\vL^2}\Bigr)^q\Bigr]\leq C_qh^{2q} + C_qk^q + C_q k^{-q}h^{2q}.
		\end{align}
		We note that the last term $C_qk^{-q}h^{2q}$ is a consequence of allowing general noises which clearly impact on the stability estimate of the time-semidiscrete pressure approximation as seen from Lemma \ref{stability_pressures}.
		
		In summary, substituting the estimates for $A_h$ and $B_h$ into \eqref{equu4.35}, we obtain
		\begin{align}
			\mE\biggl[\bigl(A_h +& B_h\bigr)^{\frac{q}{2}}\exp\Bigl(\frac{qk}{2}\sum_{n=1}^M D_h^n\Bigr)\biggr] \leq C_q\Bigl(k^{\frac{q}{2}} + h^q + k^{-\frac{q}{2}}h^q\Bigr).
		\end{align}
		The proof is complete.
	\end{proof}
	
	We conclude this section by stating a pathwise error estimate for the velocity approximation generated by Algorithm 2 which is a direct corollary of the Kolmogorov Theorem (cf. Theorem \ref{kolmogorov}) and 
	the high moment error estimates of Theorem \ref{theorem_fully_chapter5}.
	\begin{theorem}\label{theorem_pathwise_fully_chapter5}
		Assume that the assumptions of Theorem \ref{theorem_fully_chapter5} hold. Let $2\leq q <\infty$ and $0 < \gamma_2 < 1 - \frac{1}{q}$. Then, there exists a random variable $K = K(\omega)$ with $\mE\bigl[|K|^q\bigr] < \infty$ such that there holds $\mP$-a.s.
		\begin{align}
			\max_{1\leq n \leq M}\|\vu^n - \vu^n_h\|_{\vL^2} + \Bigl(\nu k \sum_{n=1}^M\|\nab(\vu^n &- \vu_h^n)\|^2_{\vL^2}\Bigr)^{1/2}  \\\nonumber
			&\leq K\bigl(k^{\frac{\gamma_2}{2}} + h^{\gamma_2} + k^{-\frac{\gamma_2}{2}}h^{\gamma_2}\bigr).
		\end{align}
	\end{theorem} 

	\subsection{High moment and pathwise error estimates for the fully discrete pressure approximation}
	In this subsection, we establish high moment and pathwise error estimates for the pressure approximation generated by Algorithm 2.  
	\begin{theorem}\label{fully_pressure_error}
		Let $2\leq q <\infty$, under the assumptions of Theorem \ref{theorem_fully_chapter5}, there holds
		\begin{align}\label{equu4.24}
			\biggl(\mE\biggl[\biggl\|k\sum_{n=1}^M\big(p^n - p^n_h\big)\biggr\|^{q}_{L^2}\biggr]\biggr)^{\frac1q} \leq C_q\bigl(k^{\frac12} + h + k^{-\frac{1}{2}}h\bigr),
		\end{align}
		where $C_q = C(\beta_1, q,D_T, \vf,\vu_0)>0$ and independent of $k$ and $h$.
	\end{theorem}

\begin{proof}
		The proof of this theorem is similar to that of Theorem \ref{semi_pressure_error}. 
		We recall the discrete inf-sup (LBB) condition: There exists a $\beta_1 > 0$ such that 
		\begin{align}
			\beta_1\|\xi_h\|_{L^2} \leq \sup_{\pphi_h\in \mH_h}\frac{\bigl(\xi_h,\div \pphi_h\bigr)}{\|\nab\pphi_h\|_{\vL^2}}\,\,\,\qquad\forall \xi_h\in L_h.
		\end{align}
		Set $\vE^m := \vu^m - \vu^m_h$ and $E_p^m := k\sum_{n = 1}^m (p^n - p^n_h)$.
		
		Summing \eqref{equu4.14} from $1$ to $M$, we obtain
		\begin{align}
			\bigl(\vE^M - \vE^0,\pphi_h\bigr) &+ \nu k\sum_{n = 1}^M\bigl(\nab\vE^n,\nab\pphi_h\bigr) + k\sum_{n=1}^M\bigl(\vu^n\cdot\nab\vu^n - \vu_h^n\cdot\nab\vu_h^n, \pphi_h\bigr) \\\nonumber
			&-\frac{k}{2}\sum_{n=1}^M\bigl([\div\vu^n_h]\vu_h^n,\pphi_h\bigr) -\bigl(E_p^M,\div \pphi_h\bigr) \\\nonumber
			&= \Bigl(\sum_{n=1}^M(\vg^{n-1} - \vg^{n-1}_h)\Delta W_{n},\pphi_h\Bigr)\,\,\qquad\forall\pphi_h\in \mH_h.
		\end{align}
		Thus,
		\begin{align}
			\bigl(E_p^M,\div \pphi_h\bigr) &= \bigl(\vE^M - \vE^0,\pphi_h\bigr) + \nu k\sum_{n = 1}^M\bigl(\nab\vE^n,\nab\pphi_h\bigr) \\\nonumber
			&\qquad+ k\sum_{n=1}^M\bigl(\vu^n\cdot\nab\vu^n - \vu_h^n\cdot\nab\vu_h^n, \pphi_h\bigr) \\\nonumber
			&\qquad-\frac{k}{2}\sum_{n=1}^M\bigl([\div\vu^n_h]\vu_h^n,\pphi_h\bigr) - \Bigl(\sum_{n=1}^M(\vg^{n-1} - \vg^{n-1}_h)\Delta W_{n},\pphi_h\Bigr)\\\nonumber
			&= \bigl(\vE^M - \vE^0,\pphi_h\bigr) + \nu k\sum_{n = 1}^M\bigl(\nab\vE^n,\nab\pphi_h\bigr) \\\nonumber
			&\qquad+ k\sum_{n=1}^M\bigl(\vE^n\cdot\nab\vu^n + \vu_h^n\cdot\nab\vE^n, \pphi_h\bigr) \\\nonumber
			&\qquad+ \frac{k}{2}\sum_{n = 1}^M \bigl([\div \vE^n]\vu_h^n,\pphi_h\bigr) - \Bigl(\sum_{n=1}^M(\vg^{n-1} - \vg^{n-1}_h)\Delta W_{n},\pphi_h\Bigr)\\\nonumber
			&\leq C(\|\vE^M\|_{\vL^2} + \|\vE^0\|_{\vL^2})\|\nab\pphi_h\|_{\vL^2} + \nu k \sum_{n = 1}^M\|\nab\vE^n\|_{\vL^2}\|\nab\pphi_h\|_{\vL^2}\\\nonumber
			&\qquad+ k\sum_{n=1}^M(\|\vE^n\|_{\vL^4}\|\nab\vu^n\|_{\vL^2} + \|\vu_h^n\|_{\vL^4}\|\nab\vE^n\|_{\vL^2})\|\pphi_h\|_{\vL^4}\\\nonumber
			&\qquad+ \frac{k}{2}\sum_{n = 1}^M\|\nab\vE^n\|_{\vL^2}\|\vu^n_h\|_{\vL^4}\|\pphi_h\|_{\vL^4}\\\nonumber
			&\qquad+ C\Bigl\|\sum_{n = 1}^M(\vg^{n-1} - \vg_h^{n-1})\Delta W_n\Bigr\|_{\vL^2}\|\nab\pphi_h\|_{\vL^2}\\\nonumber
			&\leq C(\|\vE^M\|_{\vL^2} + \|\vE^0\|_{\vL^2})\|\nab\pphi_h\|_{\vL^2} + \nu k \sum_{n = 1}^M\|\nab\vE^n\|_{\vL^2}\|\nab\pphi_h\|_{\vL^2}\\\nonumber
			&\qquad+Ck\sum_{n=1}^M(\|\nab\vE^n\|_{\vL^2}\|\nab\vu^n\|_{\vL^2} \\\nonumber
			&\qquad+ \|\nab\vu_h^n\|_{\vL^2}\|\nab\vE^n\|_{\vL^2})\|\nab\pphi_h\|_{\vL^2}\\\nonumber
			&\qquad+ \frac{Ck}{2}\sum_{n = 1}^M\|\nab\vE^n\|_{\vL^2}\|\nab\vu^n_h\|_{\vL^2}\|\nab\pphi_h\|_{\vL^2}\\\nonumber
			&\qquad+ C\Bigl\|\sum_{n = 1}^M(\vg^{n-1} - \vg_h^{n-1})\Delta W_n\Bigr\|_{\vL^2}\|\nab\pphi_h\|_{\vL^2}.
		\end{align}
		Using the discrete inf-sup condition, we obtain
		\begin{align}\label{equu4.45}
			\beta_1\|E_p^M\|_{L^2} &\leq \sup_{\pphi_h\in \mH_h}\frac{\bigl(E_p^M,\div \pphi_h\bigr)}{\|\nab\pphi_h\|_{\vL^2}}\\\nonumber
			&\leq C(\|\vE^M\|_{\vL^2} + \|\vE^0\|_{\vL^2}) + \nu k \sum_{n = 1}^M\|\nab\vE^n\|_{\vL^2}\\\nonumber
			&\qquad+Ck\sum_{n=1}^M\|\nab\vE^n\|_{\vL^2}(\|\nab\vu^n\|_{\vL^2} + \|\nab\vu_h^n\|_{\vL^2})\\\nonumber
			&\qquad+ \frac{Ck}{2}\sum_{n = 1}^M\|\nab\vE^n\|_{\vL^2}\|\nab\vu^n_h\|_{\vL^2} \\\nonumber
			&\qquad+ C\Bigl\|\sum_{n = 1}^M(\vg^{n-1} - \vg_h^{n-1})\Delta W_n\Bigr\|_{\vL^2}.
		\end{align}
		Taking the $qth$-power, expectation, and then using Theorem \ref{theorem_fully_chapter5} and the stability estimates of Lemmas \ref{stability_means} and \ref{stability_FEMs}, we obtain
		\begin{align}
			\nonumber		\beta_1\mE\bigl[\|E_p^M\|^q_{L^2}\bigr] &\leq C_q\Bigl(k^{\frac{q}{2}} + h^q + k^{-\frac{q}{2}}h^{q}\Bigr) \\\nonumber
			&\qquad+ C_q\mE\biggl[\Bigl(k\sum_{n = 1}^M\|\nab\vE^n\|^2_{\vL^2}\Bigr)^{\frac{q}{2}}\Bigl(k\sum_{n = 1}^M\bigl(\|\nab\vu^n\|^2_{\vL^2} + \|\nab\vu_h^n\|^2_{\vL^2}\bigr)\Bigr)^{\frac{q}{2}}\biggr]\\\nonumber
			&\qquad+ C_q\mE\biggl[\Bigl(k\sum_{n = 1}^M\|\nab\vE^n\|^2_{\vL^2}\Bigr)^{\frac{q}{2}}\Bigl(k\sum_{n = 1}^M\|\nab\vu_h^n\|^2_{\vL^2}\Bigr)^{\frac{q}{2}}\biggr]\\\nonumber
			&\qquad+ C_q\mE\Bigl[\Bigl(k\sum_{n = 1}^M\|\vg^{n-1} - \vg^{n-1}_h\|^2_{\vL^2}\Bigr)^{\frac{q}{2}}\Bigr]\\\nonumber
			&\leq C_q\Bigl(k^{\frac{q}{2}} + h^q + k^{-\frac{q}{2}}h^{q}\Bigr) \\
			&\qquad+ C_q\Bigl(\mE\Bigl[\Bigl(k\sum_{n = 1}^M\|\nab\vE^n\|^2_{\vL^2}\Bigr)^q\Bigr]\Bigr)^{\frac12}\Bigl(\mE\Bigl[\Bigl(k\sum_{n = 1}^M\bigl(\|\nab\vu^n\|^2_{\vL^2} + \|\nab\vu_h^n\|^2_{\vL^2}\bigr)\Bigr)^q\Bigr]\Bigr)^{\frac12}\\\nonumber
			&\qquad+ C_q\Bigl(\mE\Bigl[\Bigl(k\sum_{n = 1}^M\|\nab\vE^n\|^2_{\vL^2}\Bigr)^q\Bigr]\Bigr)^{\frac12}\Bigl(\mE\Bigl[\Bigl(k\sum_{n = 1}^M\|\nab\vu_h^n\|^2_{\vL^2}\Bigr)^q\Bigr]\Bigr)^{\frac12}\\\nonumber
			&\qquad+ C_q\mE\Bigl[\Bigl(k\sum_{n = 1}^M\|\vg^{n-1} - \vg^{n-1}_h\|^2_{\vL^2}\Bigr)^{\frac12}\Bigr]\\\nonumber
			&\leq C_q\Bigl(k^{\frac{q}{2}} + h^q + k^{-\frac{q}{2}}h^{q}\Bigr).
		\end{align}
		Finally, the proof is complete after dividing both sides by $\beta_1$.
	\end{proof}

	An immediate consequence of the above high moment error estimates is the following pathwise error estimate for the pressure approximation $\{p_h^n\}$, it follows immediately from an application of Theorem \ref{kolmogorov}.  
	
	\begin{theorem}\label{cor_4.5_chapter5}
		Assume that the assumptions of Theorem \ref{fully_pressure_error} hold. Let $2\leq q <\infty$ and $0 < \gamma_2 < 1 - \frac{1}{q}$. Then, there exists a random variable $K = K(\omega)$ with $\mE\bigl[|K|^q\bigr] < \infty$ such that there holds $\mP$-a.s.
		\begin{align}
			\Bigl\|k\sum_{n=1}^M\bigl(p^n - p^n_h\bigr)\Bigr\|_{L^2} \leq K\biggl(k^{\frac{\gamma_2}{2}} + h^{\gamma_2} + k^{-\frac{\gamma_2}{2}}h^{\gamma_2}\biggr).
		\end{align}
	\end{theorem}
	
	 \smallskip
	We conclude this section by stating the global error estimates for our fully discrete numerical 
	solution generated by Algorithm 2, which follows from combining the above temporal and spatial error estimates. 
	
	\begin{theorem} \label{thm4.8}
		Let $2 \leq q <\infty$ and $0< \gamma <\frac12$, under the assumptions of Theorem \ref{theorem_semi_chapter5} and Theorem \ref{theorem_fully_chapter5}, there exists a constant $C_q = C(D_T, \vu_0, q, \vf)>0$ such that
		\begin{align}\label{eq4.47}
			\bigl(\mE\bigl[\max_{1\leq n \leq M}\|\vu(t_n) - \vu^n_h\|^q_{\vL^2}\bigr]\bigr)^{\frac1q} &+ \biggl(\mE\biggl[\Bigl(\nu k \sum_{n=1}^{M}\|\nab\bigl(\vu(t_n) - \vu^n_h\|\bigr)\|^2_{\vL^2}\Bigr)^{\frac{q}{2}}\biggr]\biggr)^{\frac1q} \\\nonumber
			&\leq C_q\Bigl( k^{\frac12 - \gamma} + h + k^{-\frac{1}{2}}h\Bigr).
		\end{align}
		In addition, let $2<q<\infty$ and $0<\gamma < \frac12$ such that $\frac12 - \gamma - \frac{1}{q}>0$. Then, for any $0<\gamma_1 < \frac12 - \gamma - \frac{1}{q}$ and $0<\gamma_2 < 1 - \frac{1}{q}$, there exists a random variable $K$ with $\mE[|K|^q] < \infty$ such that there holds $\mP-a.s.$
		\begin{align}\label{equu5.2}
			\max_{1\leq n \leq M}\|\vu(t_n) - \vu^n_h\|_{\vL^2} &+ \Bigl(\nu k \sum_{n=1}^{M}\|\nab\bigl(\vu(t_n) - \vu^n_h\bigr)\|^2_{\vL^2}\Bigr)^{\frac12} \\\nonumber
			&\leq K\Bigl(k^{\gamma_1} + h^{\gamma_2} + k^{-\frac{\gamma_2}{2}}h^{\gamma_2}\Bigr). 
		\end{align}
	\end{theorem}

	\begin{theorem} \label{thm4.9}
		Let $2 \leq q <\infty$ and $0< \gamma <\frac12$. Under the assumptions of Theorem \ref{semi_pressure_error} and Theorem \ref{fully_pressure_error}, there exists a constant $C_q = C(D_T, \vu_0, q, \vf,\beta_0,\beta_1)>0$ such that for $1 \leq \ell \leq M$ 
		\begin{align}\label{equu4.33}
			\biggl(\mE\biggl[\biggl\|P(t_{\ell}) - k\sum_{n=1}^{\ell} p^n_h\biggr\|^q_{L^2}\biggr]\biggr)^{\frac1q} \leq C_q\Bigl( k^{\frac12 - \gamma} + h + k^{-\frac{1}{2}}h\Bigr),
		\end{align}
		In addition, let $2<q<\infty$ and $0<\gamma < \frac12$ such that $\frac12 - \gamma - \frac{1}{q}>0$. Then, for any $0<\gamma_1 < \frac12 - \gamma - \frac{1}{q}$ and $0<\gamma_2 < 1 - \frac{1}{q}$, there exists a random variable $K$ with $\mE[|K|^q] < \infty$ such that there holds $\mP-a.s.$
		\begin{align}\label{eq5.104}
			\biggl\|P(t_{\ell}) - k\sum_{n=1}^{\ell} p^n_h\biggr\|_{L^2} \leq K\Bigl(k^{\gamma_1} + h^{\gamma_2} + k^{-\frac{\gamma_2}{2}}h^{\gamma_2}\Bigr).
		\end{align}
	\end{theorem}
	
	\section{Numerical experiments} \label{sec-5}
		In this section, we present three numerical experiments to validate our error estimates in \eqref{eq4.47}--\eqref{eq5.104}. 
		In all our experiments,  we set $D = (0,1)^2\subset \mathbb{R}^2$, $T =1, \nu =1$, $g \equiv 10$, the body force is $\vf = (f_1,f_2)$ with
		\begin{align*}
			f_1(x,y) &= \pi\cos(t)\sin(2\pi y)\sin(\pi x)\sin(\pi x)
			- 2\pi^3\sin(t)\sin(2\pi y)(2\cos(2\pi x)-1) \\\nonumber
			&\qquad- \pi\sin(t)\sin(\pi x)\sin(\pi y),\\\nonumber
			f_2(x,y) &= -\pi\cos(t)\sin(2\pi x)\sin(\pi y)\sin(\pi y)
			- 2\pi^3\sin(t)\sin(2\pi x)(1-2\cos(2\pi y)) \\\nonumber
			&\qquad+ \pi\sin(t)\cos(\pi x)\cos(\pi y).
		\end{align*}
		We choose $\vW$ in \eqref{eq1.1} to be a ${\mathbb R}^J$-valued Wiener process, with  increment
		\begin{align}\label{used_noise}
			\vW^J(t_{n+1},{\bf x}) - \vW^J(t_n,{\bf x}) = \sqrt{k_0} \sum\limits_{j_1=1}^J\sum\limits_{j_1=1}^J \sqrt{\lambda_{j_1,j_2}}{\bf e}_{j_1,j_2}({\bf x})\xi_{j_1,j_2}^n\,,
		\end{align}
		where ${\bf x}=(x,y) \in D, \, \xi_{j_1,j_2}^n \sim N(0,1)$, $\lambda_{j_1,j_2} = \frac{1}{(j_1+j_2)^2}$, and
		\begin{align}
			{\bf e}_{j_1,j_2}({\bf x}) = \bigl(\sin(j_1\pi x)\sin(j_2 \pi y),\sin(j_1\pi x)\sin(j_2 \pi y)\bigr)\, .
		\end{align}
		In addition, we use the Taylor-Hood mixed finite element method for the spatial discretization and the homogeneous Dirichlet boundary condition is imposed on $\vu$. 
		Moreover, we choose $J = 4$ and $N_p =300$ as number of samples in our Monte Carlo simulations to compute the expectations in each test. 
		We implement Algorithm 2 and evaluate the errors between the computed velocity (resp. pressure) and the reference velocity (resp. pressure) which are generated by Algorithm 2 on a very fine mesh with mesh sizes $k_0 = \frac{1}{1024}, h =\frac{1}{40}$. Furthermore, to evaluate errors in strong norms, we use the following numerical 
		integration formulas: for any $0 \leq m \leq M$ and any integer $2\leq q <\infty$,
		\begin{align*}
			\pmb{\mathcal{E}}_{\vu,q}^m &:=\biggl(\mE\Bigl[\|{\bf u}(t_m) -{\bf u}_h^m(k)\|_{\vL^2}^q \Bigr]\biggr)^{\frac1q} \approx
			\biggl(\frac{1}{N_p}\sum_{\ell=1}^{N_p} \bigl\|{\bf u }_h^m(k_0,\omega_\ell)
			-{\bf u}_h^m(k,\omega_\ell) \bigr\|_{\vL^2}^q \biggr)^{\frac1q}\, ,\\\nonumber
			{\mathcal{E}}_{P,q}^m &:=\biggl(\mE\Bigl[\Big\|P(t_m) -k\sum_{n=1}^m p_h^n(k)\Big\|_{\vL^2}^q \Bigr]\biggr)^{\frac1q}
			\\\nonumber
			&\,\,\approx \biggl(\frac{1}{N_p}\sum_{\ell=1}^{N_p}\Bigl\|k\sum_{n=1}^{\frac{t_m}{k_0}} p_{h}^{n}(k_0,\omega_\ell) -k\sum_{n=1}^{\frac{t_m}{k}}p_h^n(k,\omega_\ell)\Bigr\|_{\vL^2}^q\Bigr) \biggr)^{\frac1q}\, .
		\end{align*}
		
		In addition, in order to compute the quantities of stochastic interest, we want to generate as many samples as possible in the Monte Carlo simulations. However, this approach is very expensive if Newton's method is used at each time step to solve the nonlinear system in Algorithm 2. Therefore, we introduce a cheaper fixed point iteration given below  in Algorithm 3  to solve the nonlinear system.
		
		\medskip 
		{\bf Algorithm 3.} Given $\vu_h^{n+1,\ell-1} = \vu_h^n$ for $n = 0, 1,\cdots, M-1$; $\ell = 1,\cdots,L$. Find $(\vu_h^{n+1,\ell}, p_h^{n+1,\ell}) \in L^2(\Ome;\mH_h\times L_h)$ such that $\mP$-a.s.
		\begin{align*}
			\bigl(\vu_h^{n+1,\ell},\pphi_h\bigr) &+ \nu k\bigl(\nab\vu_h^{n+1,\ell},\pphi_h\bigr) - k\bigl(p_h^{n+1,\ell},\div\pphi_h\bigr) \\\nonumber
			&= \bigl(\vu_h^{n},\pphi_h\bigr) - k\tilde{b}\bigl(\vu_h^{n+1,\ell-1},\vu_h^{n+1,\ell-1},\pphi_h\bigr) + \bigl(g_h^n\Delta \vW_{n+1},\pphi_h\bigr)\,\,\forall \pphi_h\in\mH_h,\\
			\bigl(\div\vu_h^{n+1,\ell},\psi_h\bigr) &=0\qquad\forall \psi_h\in L_h.
		\end{align*}
		
		\medskip
		{\bf Test 1.} In this test, we want to verify the convergence rates in the error estimates \eqref{eq4.47}, \eqref{equu4.33}. To do that, the high moment errors $\E_{\vu,q}^M$ and $\E_{P,q}^M$ with $q = 2,4,8$ of the velocity and pressure approximations generated by Algorithm 2 are computed and the numerical results are shown in Table \ref{table5.1} and  \ref{table5.2}. The numerical results verify almost a half order convergence rate as predicted by our error estimate results. 
		\begin{table}[tbhp]
			\begin{center}
				\begin{tabular}{ |c|c|c|c|c|c|c|}
					\hline
					\bf $k$ & $\pmb{\mathcal{E}}_{\vu,2}^M$  & \mbox{order} & $\pmb{\mathcal{E}}_{\vu,4}^M$ &  \mbox{order}&$\pmb{\mathcal{E}}_{\vu,8}^M$&order\\
					\hline 
					$1/64$  & 0.01031 &  & 0.0136263&&0.0189242&\\
					\hline
					$1/128$  & 0.00873561 &0.2357  & 0.0118123&0.2068&0.0169163&0.1637\\
					\hline
					$1/256$  & 0.00618757 & 0.4975 & 0.007929&0.5751&0.010652&0.6673\\
					\hline
					$1/512$  & 0.0045004 & 0.4594 & 0.00576594&0.4586&0.00776608&0.4555\\
					\hline
				\end{tabular}
				\caption{Algorithm 2: Errors of the computed velocity $\{ {\bf u}^n_h\}_n$.}
				\label{table5.1}
			\end{center}
		\end{table}
		\begin{table}[tbhp]
			\begin{center}
				\begin{tabular}{ |c|c|c|c|c|c|c|}
					\hline
					\bf $k$ & ${\mathcal{E}}_{P,2}^M$  & \mbox{order} & ${\mathcal{E}}_{P,4}^M$ &  \mbox{order}&${\mathcal{E}}_{P,8}^M$&order\\
					\hline 
					$1/64$  & 0.048219 &  & 0.0565676&&0.0648501&\\
					\hline
					$1/128$  & 0.0335814 & 0.5219 & 0.039123&0.5320&0.0446196&0.5397\\
					\hline
					$1/256$  & 0.02364 & 0.5064 & 0.0274858&0.5093&0.0312931&0.5118\\
					\hline
					$1/512$  &  0.0155695& 0.6025 &0.0180596 &0.6060&0.020537&0.5758\\
					\hline
				\end{tabular}
				\caption{Algorithm 2: Errors of the computed pressure $\{ {p}^n_h\}_n$.}
				\label{table5.2}
			\end{center}
		\end{table}
		
\begin{figure}[thb]
			\begin{center}
				\centerline{
					\includegraphics[scale=0.4]{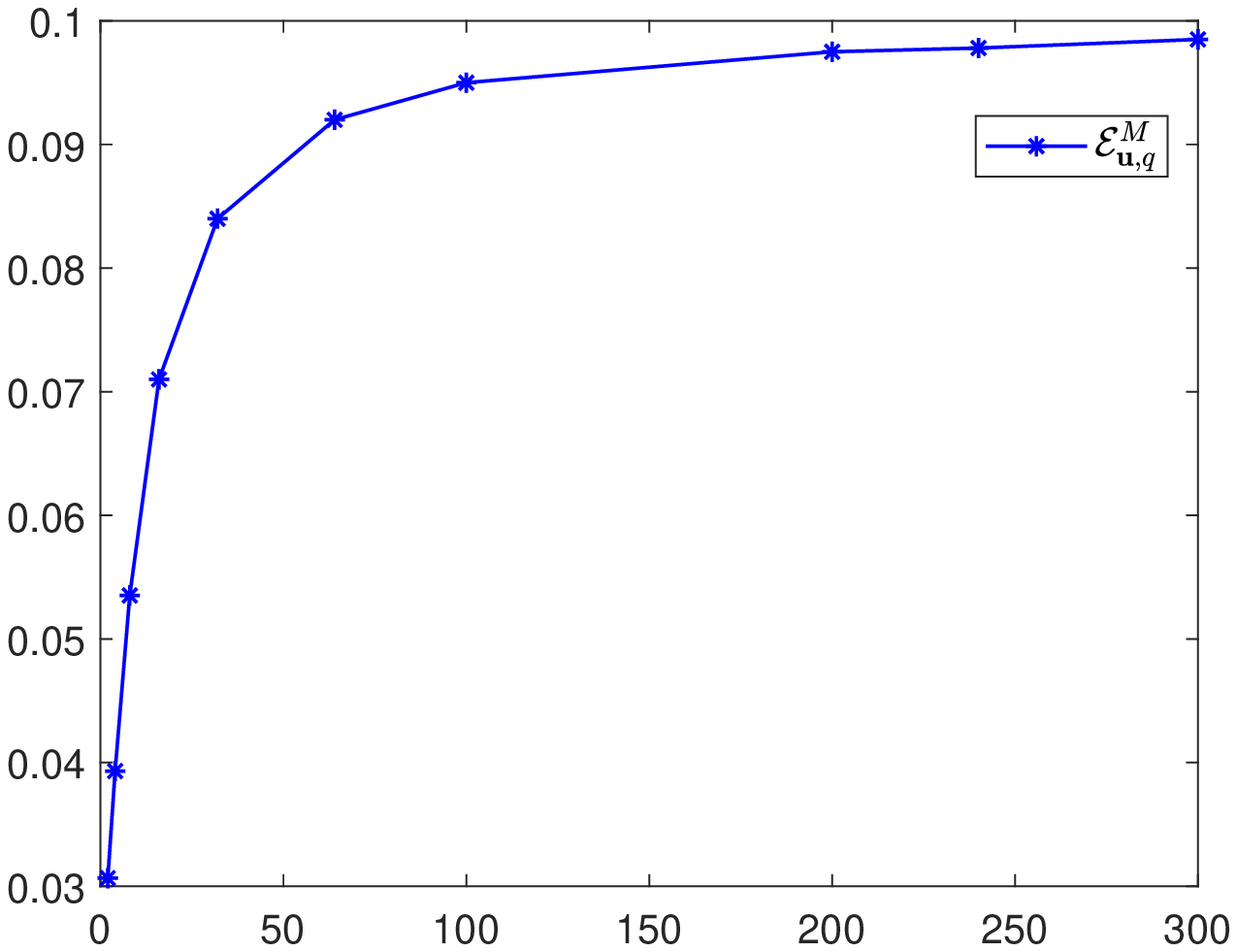}
					\includegraphics[scale=0.4]{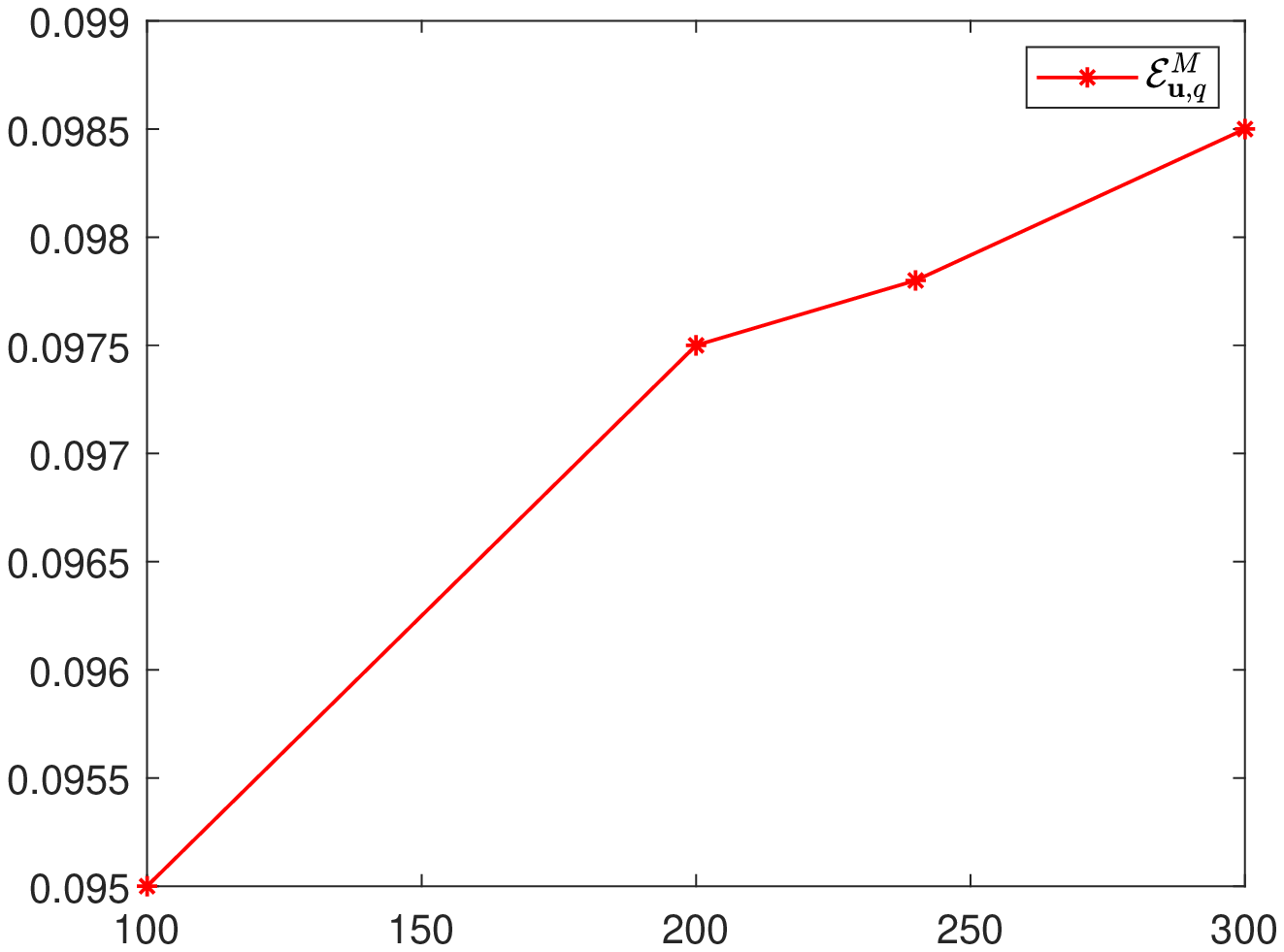}
				}
				\caption{Errors of the velocity approximation (left) in $\pmb{\E}^M_{\vu,q}$ norm and its behavior (right) for large $q$.}\label{figure5.1}
			\end{center}
		\end{figure}

\begin{figure}[thb]
			\begin{center}
				\centerline{
					\includegraphics[scale=0.4]{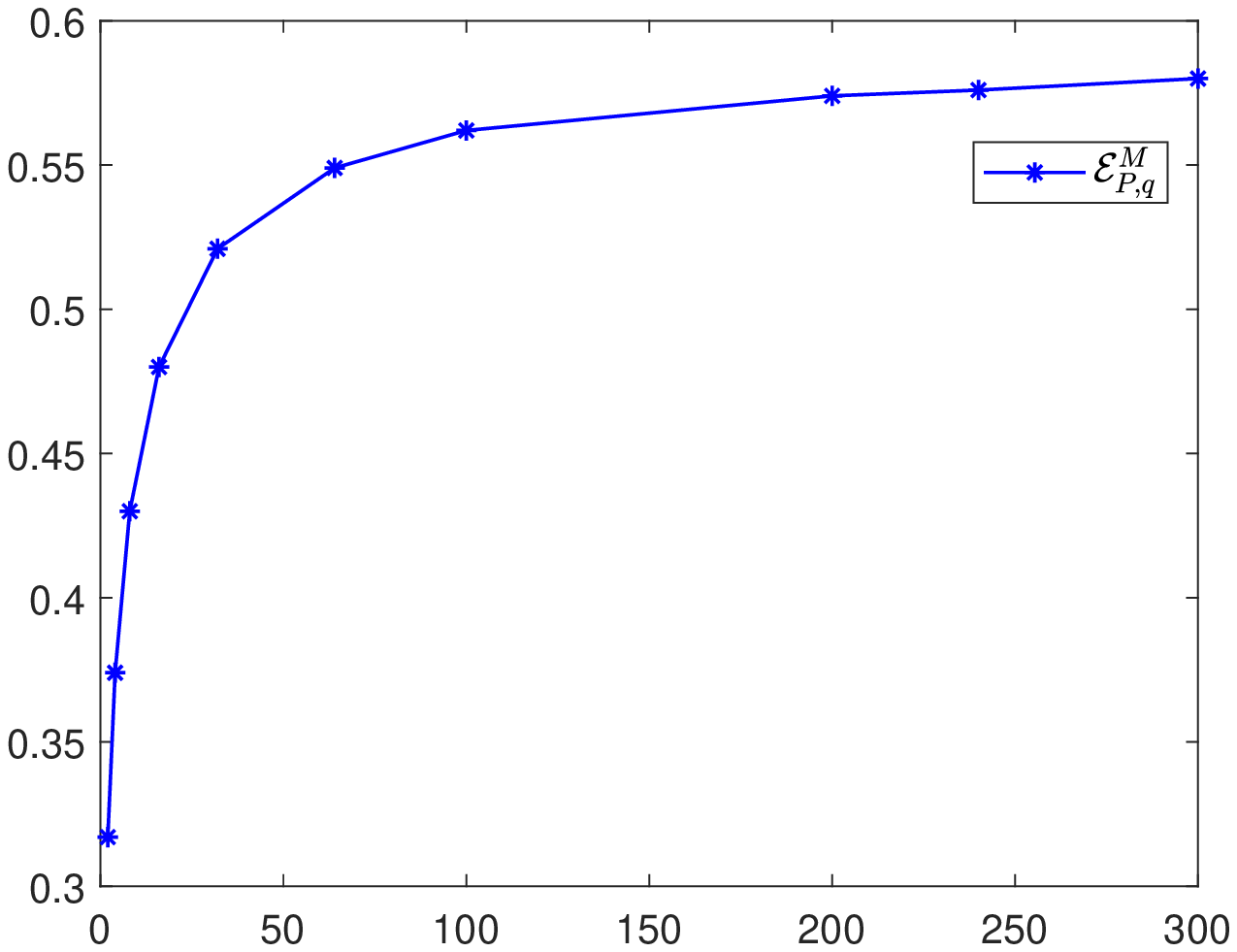}
					\includegraphics[scale=0.4]{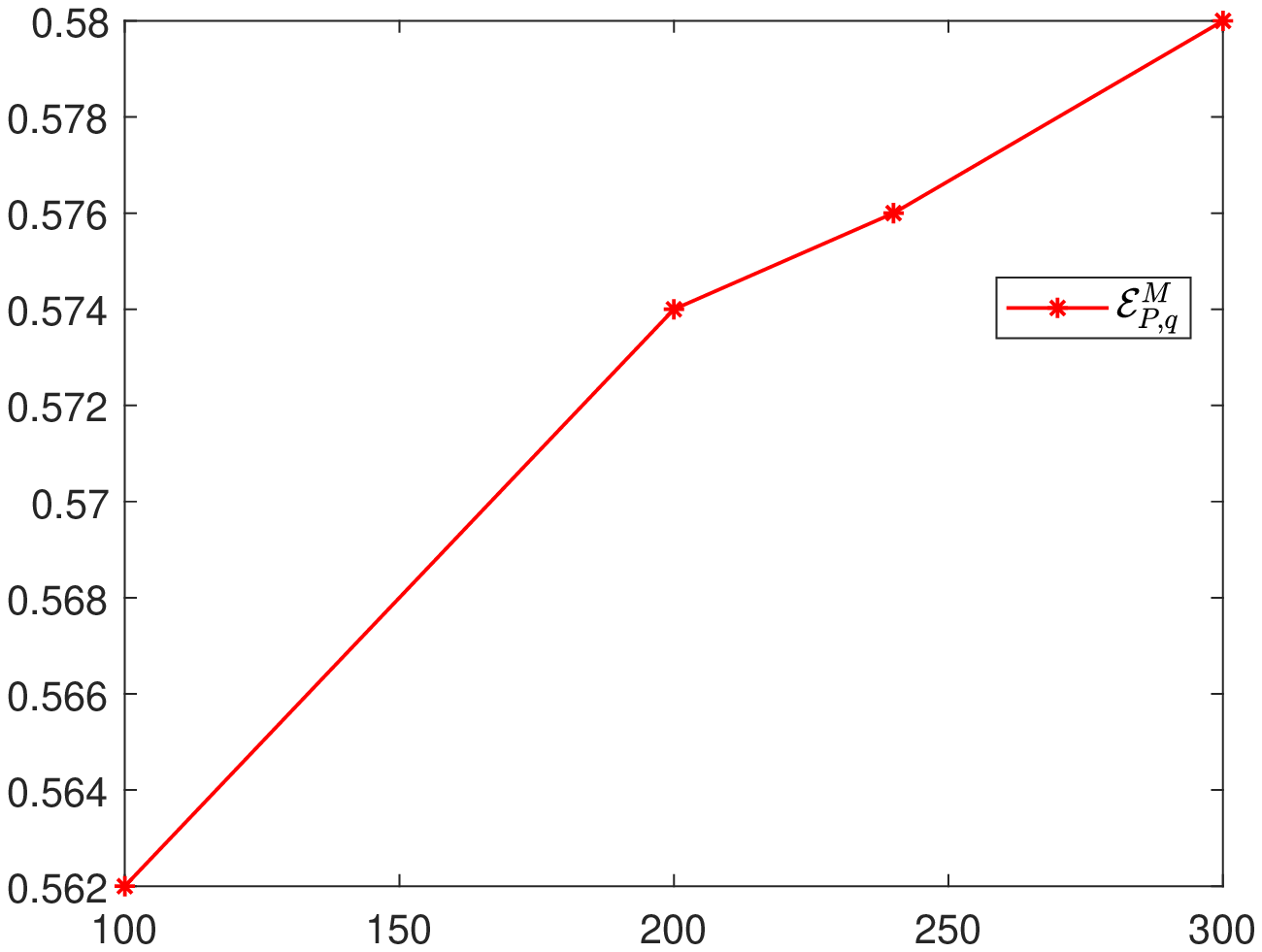}
				}
				\caption{Errors of the pressure approximation (left) in $\E_{P,q}^M$ norm and its  behavior (right)  for large $q$.}\label{figure5.2}
			\end{center}
		\end{figure}

		{\bf Test 2.} The constant $C_q$ in Theorem \ref{thm4.8} and Theorem \ref{thm4.9} depends on $q$. However, no analytical formula for $C_q$ is known. The goal of this test is to find out numerically how the error constant $C_q$ depends on $q$. To the end, we fix $h = \frac{1}{40}$ and choose the time step $k = \frac{1}{16}$ to compute the errors $\E_{\vu,q}^M$ and $\E_{P,q}^M$ for different values of $q$. The numerical results are given  in Figure \ref{figure5.1} and \ref{figure5.2} (left). We observe that the numerical results  suggest the constant $C_q$ is increasing (and blows up) in $q$. For a close-up view, we zoom in at large
		$q$ to see the behaviors of both velocity and pressure errors in Figure \ref{figure5.1} and \ref{figure5.2} (right). We still see the increasing  of the errors in $q$ although the  growth becomes slower for large $q$. A consequence of this analysis also shows that we can not simply take limit as $p\to \infty$ in the high moment error estimates of Theorem \ref{thm4.8} and \ref{thm4.9} to derive pathwise error estimates \eqref{equu5.2} and \eqref{eq5.104}, and using Kolmogorov's Theorem is still the only viable approach for obtaining pathwise error estimates.  
		
\begin{figure}[thb]
			\begin{center}
				\centerline{
					\includegraphics[scale=0.4]{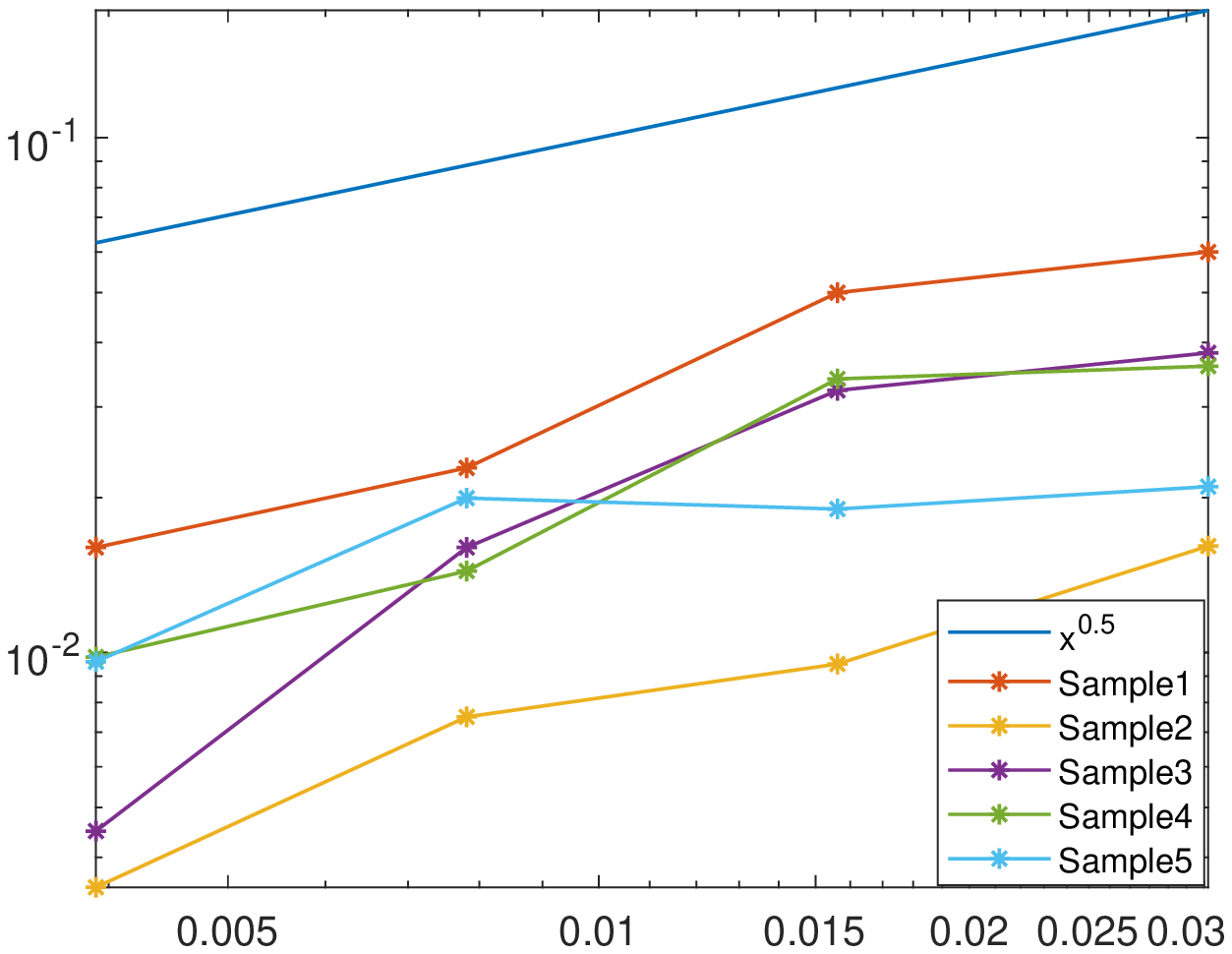}
					\includegraphics[scale=0.4]{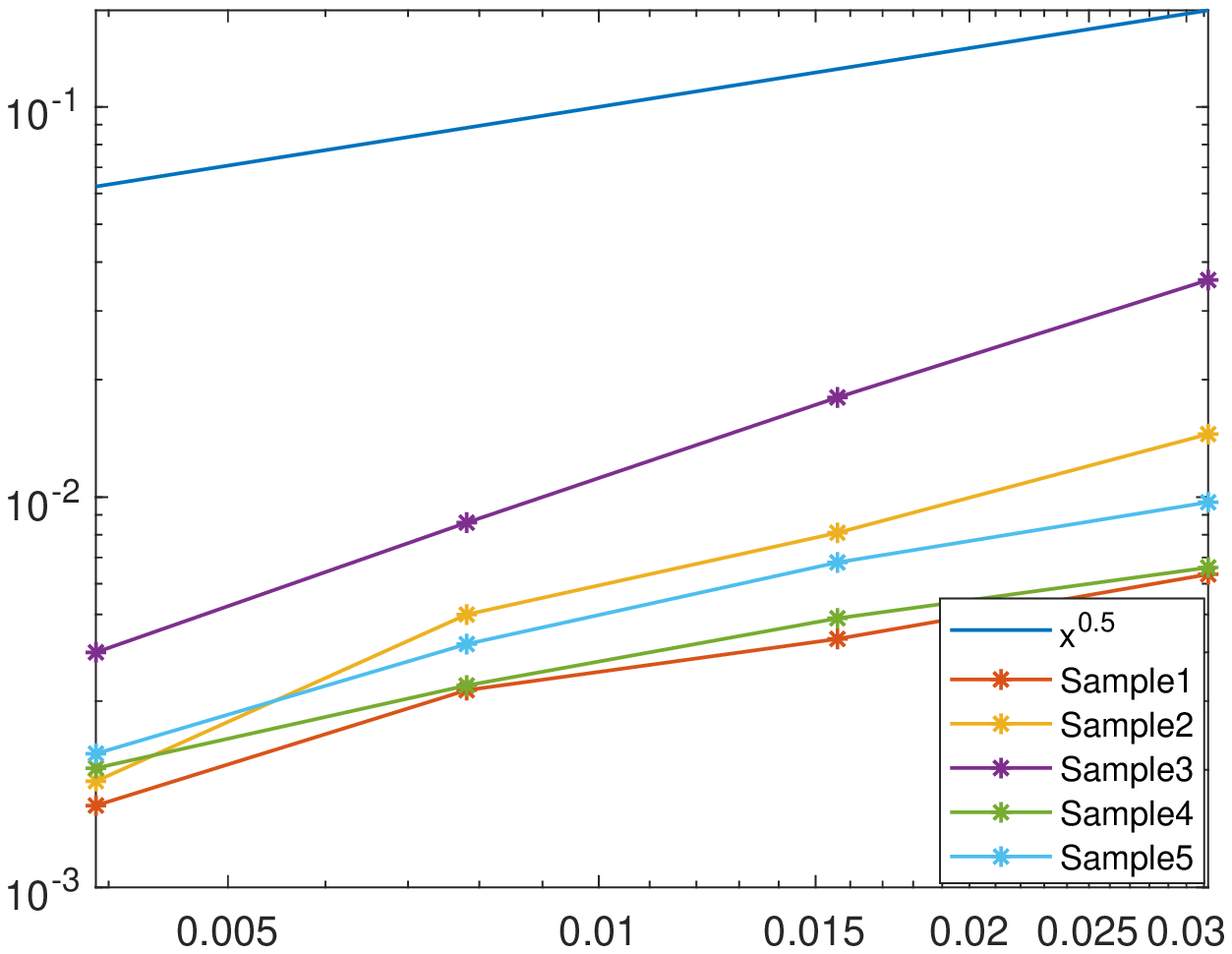}
				}
				\caption{Five sample pathwise errors of the velocity (left) and pressure approximations (right).}\label{fig5.3}
			\end{center}
		\end{figure}

		\medskip
		{\bf Test 3.} In this test, we verify the pathwise error estimates \eqref{equu5.2} and \eqref{eq5.104}. We choose five sample paths and compute their $L^2$-norm errors for both velocity and pressure approximations. The numerical results are displayed in Figure \ref{fig5.3} which verify that the pathwise convergence of the velocity and pressure approximations is almost of order $O(k^{\frac12})$ as predicted by our pathwise error estimates. 



\end{document}